\documentclass[a4paper,DIV=12,11pt,parskip=half]{scrartcl}

\usepackage{graphicx}        
\usepackage{multicol}        

\usepackage[utf8x]{inputenc}

\makeatletter
\let\@@citation@@=\citation

\renewcommand{\citation}[1]{\@@citation@@{#1}%
	\@for\@tempa:=#1\do{\@ifundefined{cit@\@tempa}%
		{\global\@namedef{cit@\@tempa}{}}{}}%
}
\makeatother

\usepackage{hyperref}

\usepackage{amssymb,amsmath}

\usepackage{mathtools}
\mathtoolsset{centercolon}

\usepackage{booktabs}

\usepackage{siunitx}

\usepackage{caption}
\usepackage{subcaption}
\usepackage{cleveref}

\usepackage{autonum}

\usepackage{enumitem}

\usepackage[section]{placeins}

\makeatletter
\def\@lbibitem[#1]#2#3\par{%
	\@ifundefined{cit@#2}{}{\@skiphyperreftrue
		\H@item[%
		\ifx\Hy@raisedlink\@empty
		\hyper@anchorstart{cite.#2\@extra@b@citeb}%
		\@BIBLABEL{#1}%
		\hyper@anchorend
		\else
		\Hy@raisedlink{%
			\hyper@anchorstart{cite.#2\@extra@b@citeb}\hyper@anchorend
		}%
		\@BIBLABEL{#1}%
		\fi
		\hfill
		]%
		\@skiphyperreffalse}%
	\if@filesw
	\begingroup
	\let\protect\noexpand
	\immediate\write\@auxout{%
		\string\bibcite{#2}{#1}%
	}%
	\endgroup
	\fi
	\ignorespaces
	\@ifundefined{cit@#2}{}{#3}}

\def\@bibitem#1#2\par{%
	\@ifundefined{cit@#1}{}{\@skiphyperreftrue\H@item\@skiphyperreffalse
		\Hy@raisedlink{%
			\hyper@anchorstart{cite.#1\@extra@b@citeb}\relax\hyper@anchorend
	}}%
	\if@filesw
	\begingroup
	\let\protect\noexpand
	\immediate\write\@auxout{%
		\string\bibcite{#1}{\the\value{\@listctr}}%
	}%
	\endgroup
	\fi
	\ignorespaces
	\@ifundefined{cit@#1}{}{#2}}
\makeatother

\renewcommand{\vec}{\boldsymbol}
\newcommand{\mat}{\boldsymbol}

\renewcommand{\d}{\mathrm d}
\renewcommand{\P}{\mathbb{P}}
\newcommand{\R}{\mathbb{R}}
\newcommand{\N}{\mathbb{N}}

\newcommand{\force}{{\mat f}}
\newcommand{\forceC}{f}

\newcommand{\displacement}{\mat u}
\newcommand{\displacementC}{u}
\newcommand{\velocity}{\mat v}
\newcommand{\velocityC}{v}

\newcommand\Tau{\mathcal{T}} 

\renewcommand{\j}{\iota}

\newcommand{\overbar}[1]{\mkern 1.5mu\overline{\mkern-2.5mu#1\mkern-0.0mu}\mkern 1.5mu}

\DeclareMathOperator{\spann}{span}
\newcommand\spanset[1]{\ensuremath\spann\{#1\}}

\newtheorem{remark}{Remark}
\newtheorem{definition}{Definition}
\newtheorem{theorem}{Theorem}

\begin{document}

\title{Numerical study of Galerkin--collocation approximation in time for the wave 
equation}
\author{Mathias Anselmann$^{\star}$, Markus Bause$^{\star}$\\
{\small ${}^\star$ Helmut Schmidt University, Faculty of 
	Mechanical Engineering, Holstenhofweg 85,}\\ 
{\small 22043 Hamburg, Germany}\\
}

\date{}

\maketitle

\begin{abstract}
\textbf{Abstract.} The elucidation of many physical problems in science and engineering is subject 
to the accurate numerical modelling of complex wave propagation phenomena. Over the last 
decades, high-order numerical approximation for partial differential equations has become 
a well-established tool. Here we propose and study numerically the implicit 
approximation in time of wave equations by a Galerkin--collocation approach that relies 
on a higher order space-time finite element approach. The conceptual basis is the 
establishment of a direct connection between the Galerkin method for the time 
discretization and the classical collocation methods, with the perspective of achieving 
the accuracy of the former with reduced computational costs provided by the latter in 
terms of less complex linear algebraic systems. For the fully discrete solution, higher 
order regularity in time is further ensured which can be advantageous in the 
discretization of multi-physics systems. The accuracy and efficiency of the variational 
collocation approach is carefully studied by numerical experiments.
\end{abstract}

\section{Introduction}

\begin{figure}[h!tb]
	\centering
	\includegraphics[width=.9\textwidth]{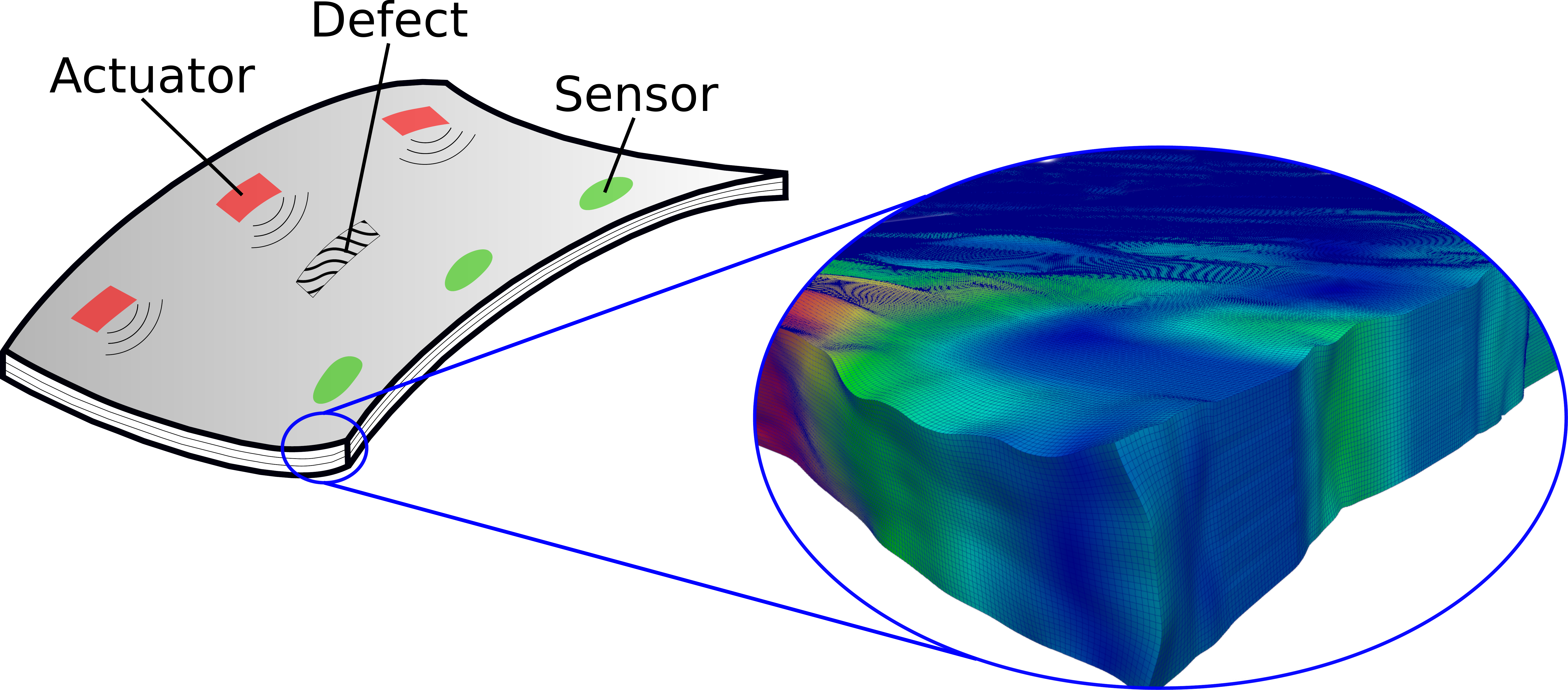}
	\caption{Concept of structural health monitoring with finite 
element simulation (scaled displacement field) illustrating the expansion of elastic 
waves.}
	\label{fig:elastic_wave}
\end{figure}

The accurate and efficient numerical simulation of wave phenomena continues to remain a 
challenging task and attract researchers' interest. Wave phenomena are studied in various 
branches of natural sciences and technology. For instance, fluid-structure interaction, 
acoustics, poroelasticity, seismics, electro-magnetics and non-destructive material 
inspection represent prominent fields in that wave propagation is studied. One of our key 
application for wave propagation is structural health monitoring of lightweight 
material (for instance, carbon-fibre reinforced polymers) by ultrasonic waves in 
aerospace engineering. The conceptional idea of this new and intelligent approach is 
sketched in \cref{fig:elastic_wave}.
The structure is equipped with an integrated actuator-sensor network. The ultrasonic 
waves that are emitted by the actuators interact with material defects of the solid 
structure. By means of an inverse modelling, the signals that are recorded by the sensors 
monitor material failure (cf. \cite{Koecher15}) and, as perspective for the future, may 
allow prognoses about the structure's residual lifetime. The design of such monitoring 
systems and the signal interpretation require the elucidation of wave propagation in 
composite material which demands for highly advanced and efficient numerical simulation 
techniques. 

High-order numerical approximation of partial differential equations has been strongly 
focused and investigated in the last decades. High order methods are known to be 
efficient if they approximate functions with large elements of high polynomial degree in 
regions of high regularity. Prominent examples are $hp$- and spectral element methods in 
application areas such as computational fluid dynamics or computational mechanics. 
Their theoretical convergence analysis and adaptive $hp$- and spectral 
element versions still experience strong development. Whereas high-order approaches have 
been considered for the approximation of the spatial variables, first- or second-order 
implicit schemes are often still used for the discretization of the time variable. We 
note that in this work only implicit time discretization schemes are in the scope of 
interest. Thus, all remarks refer to this class of methods. Our motivation for 
using implicit time discretization schemes comes from the overall goal to apply the 
proposed Galerkin--collocation techniques to mixed systems like, for 
instance, fluid-structure interaction for free flow modelled by the Navier--Stokes 
equations \cite{Richter2017} or fully dynamic poroelasticity \cite{Wheeler2012}. 

Driven by the tremendous increase in computing power of modern high performance computing 
systems and recent progress in the technology of algebraic solver, including efficient 
techniques of preconditioning, space-time finite element approaches have recently 
attracted high attention and have been brought to application maturity; cf., e.g., 
\cite{Hussain2014,Doerfler2016,Ernesti2018,Steinbach2015}. Space-time finite element methods offer 
appreciable advantages over discretizations of mixed type based on finite difference 
techniques for the discretization of the time variable (e.g., by Runge-Kutta methods) 
and, for instance, finite element methods for the discretization of the space variables. 
In particular, advantages are the natural embedding of higher order members in the 
various families of schemes, the applicability of functional analysis techniques in their 
analyses due to the uniform space-time framework and the applicability of well-known 
adaptive mesh refinement techniques, including goal-oriented error control 
\cite{Bangerth2010}. In the meantime a broad variety of implementations of space-time 
finite element methods does exist. The families of schemes 
differ by the choices of the trial and test spaces. This leads to continuous or 
discontinuous approximations of the time variable (cf. e.g., 
\cite{Matthies2011,Zhao2014}). Further, the fully coupled treatment of all time 
steps versus time-marching approaches is discussed. In particular, the simultaneous 
computation of all time steps imposes high demands on the linear solver technology (cf., 
e.g., \cite{Doerfler2016,Ernesti2017,Ernesti2018}).

In this work, we propose the Galerkin--collocation method for the numerical solution of 
wave equations. This approach combines variational approximation in time by finite 
element techniques with the concepts of collocation methods and follows the ideas of 
\cite{Becher2018}. By imposing collocation conditions, the test space of the variational 
condition is downsized. The key ingredients and innovations of the approach are:
\begin{itemize}
	\item[A.] Higher order regularity in time of the fully discrete 
approximation;  
	\item[B.] Linear systems of reduced complexity;   
\end{itemize}	
Ingredient [A] is a direct consequence of the construction of the schemes. Higher order 
regularity is ensured by imposing collocation conditions at the discrete time nodes and 
endpoints of the subintervals $[t_{n-1},t_n]$, for $n=1,\ldots , N$, of the global time 
interval $[0,T]$. Higher order regularity in time might offer appreciable advantages 
for future approximations of coupled multi-physics systems if higher order time 
derivatives of the discrete solution of one subproblem arise as coefficient functions in 
other subproblems. Ingredient  [B] is ensured by the proper choice of a special for the 
discrete in time function spaces. Thereby, simple vector identities for the degrees of 
freedom in time are obtained at the left endpoints of the subintervals without generating 
computational costs. These vector identities can then be exploited to eliminate 
conditions from the algebraic systems and reduce its size compared to the standard 
continuous Galerkin--Petrov approximation in time; cf.~\cite{Bause2018}. In a further 
work of the authors \cite{Bause2019}, it is shown that the optimal order of convergence 
in time (and space) of the underlying finite element discretization is preserved by the 
Galerkin--collocation approach. The numerical example 
of \cref{sec:Challenging_Example} that mimics typical studies of 
structural health monitoring (cf.\ \cref{fig:elastic_wave}), demonstrates the 
superiority of the Galerkin--collocation approach over a standard continuous 
Galerkin--Petrov method admitting continuity and no 
differentiability in time of the discrete solution.

For the sake of brevity, standard conforming finite element methods are used for the 
discretization of the spatial variables in this work. This is done since we focus here on 
time discretization. In the literature it has been mentioned that discontinuous 
finite element methods in space offer appreciable advantages over continuous ones for the  
discretization of wave equations; cf., e.g., \cite{Basabe2008,Grote2006}. The 
application of, for instance, the symmetric interior penalty disccontinuous Galerkin 
method (cf.\ \cite{Bause2014,Koecher15} along with a Galerkin--collocation 
discretization in time, is straightforward. 

This work is organized as follows. In \cref{sec:ProblemDescription} we 
introduce our prototype model. In \cref{sec:Gal-col_scheme} we present its 
discretization by two families of Galerkin--collocation methods.  In \cref{sec:c1_solution}, the discrete form of a member of theses families with 
$C^1$-regularity in time is derived. The resulting algebraic system is built and our 
algebraic solver is described. In \cref{sec:C2_solution} the discrete form of a member 
of the Galerkin--collocation family with $C^2$-regularity in time is derived.
For both methods, the results of our numerical experiments are presented and evaluated. 

\section{Mathematical problem and notation}
\label{sec:ProblemDescription}

As a prototype model, we study the wave problem 
\begin{equation}
\label{eq:ScalarWaveSO}
\begin{array}{r@{\;}c@{\;}ll}
\partial_t^2 u - c^2 \Delta{u}  & = f\,, & \in \Omega \times I\,,\\[1ex]
u(0) = u_0\,, \quad \partial_t u(0)  & = v_0 \,, &  \text{in}\; \Omega\,,\\[1ex]
u  = g^u \,, \;\;  \text{on}\; \partial \Omega_D \times I\,,\quad \partial_n u  & = 0 
\,, & \text{on}\; \partial \Omega_N \times I\,.\\
\end{array}
\end{equation}
In our application of structural health monitoring (cf.\ \cref{fig:elastic_wave}), $u$ denotes the scalar valued displacement field, $c\in 
R$ with $c>0$, is a material parameter and $f$ an external force acting 
on the domain 
$\Omega\subset \R^d$, with $d=2,3$. Further, $g^u$ is a prescribed trace on the Dirichlet 
part $\partial \Omega_D$ of the boundary $\partial \Omega=\partial \Omega_D \cup \partial 
\Omega_N$, with $\partial \Omega_D \cap \partial \Omega_N = \emptyset$. By $\partial_n$ 
we denote the normal derivative with outer unit normal vector $\mat n$. Homogeneous 
Neumann boundary conditions on $\partial \Omega_N$ are prescribed for brevity. Finally, 
$I=(0,T]$ denotes the time domain. Problem \eqref{eq:ScalarWaveSO} is 
well-posed and admits a unique solution $(u,\partial_t u) \in L^2(0,T;H^1(\Omega))\times 
L^2(0,T;L^2(\Omega))$ under appropriate assumptions about the data; cf.\ \cite{Lions1968}. 
By imbedding, $u \in C([0,T];H^1(\Omega))$ and $v \in C([0,T];L^2(\Omega))$ is ensured; 
cf.\ \cite{Lions1971}. Throughout, we tacitly assume that the solution admits all the 
(improved) regularity being necessary in the arguments.  

Our notation is standard. By $H^m(\Omega)$ we denote the Sobolev space of $L^2(\Omega)$ 
functions with derivatives up to order $m$ in $L^2(\Omega)$. For brevity, we let $H := 
L^2(\Omega)$ and $V=H^1_{0,D}(\Omega)$ be the space of all $H^1$-functions with 
vanishing trace on the Dirichlet part $\partial \Omega_D$ of $\partial 
\Omega$. By $\langle \cdot,\cdot \rangle_{\Omega}$ we denote the inner 
product in $L^2(\Omega)$. For the norms we use $\| \cdot \| := \| \cdot\|_{L^2(\Omega)}$ 
and $\| \cdot \|_m := \| \cdot\|_{H^m(\Omega)}$ for $m\in \N$ and $m\geq 1$. Finally, the 
expression $a\lesssim b$ stands for the inequality $a \leq C\, b$ with a generic constant 
$C$ that is indepedent of the size of the space and time meshes. 

By $L^2(0,T;B)$, $C([0,T];B)$ and $C^q([0,T];B)$, for $q\in \N$, we denote the 
standard Bochner spaces of $B$-valued functions for a Banach space $B$, equipped with 
their natural norms. Further, for a subinterval $J\subseteq [0,T]$, we will use the 
notations $L^2(J;B)$, $C^m(J;B)$ and $C^0(J;B):= C(J;B)$ for the corresponding Bochner 
spaces.

To derive our Galerkin--collocation approach, we first rewrite problem 
\eqref{eq:ScalarWaveSO} as a first order system in time for the unknowns $(u,v)$, with 
$v=\partial_t u$,
\begin{equation}
\label{eq:InitialProblem}
\partial_t u - v  = 0\,, \quad \partial_t v - c^2  \Delta{u}  = f\,.
\end{equation}
Further, we represent the unknowns $u$ and $v$ in terms of 
\begin{equation}
 \label{eq:SplitSol}
 u = u^0 + u^D \quad \text{and} \quad v = v^0 + v^D\,.
\end{equation}
Here, $u^D$, $v^D\in C(\overline I;H^1(\Omega))$ are supposed to be (extended) functions 
with traces $u^D=g^u$ and $v^D=g^v:=\partial_t g^u$ on the Dirichlet part $\partial 
\Omega_D$ of 
$\partial \Omega$.    

Using \eqref{eq:InitialProblem} and \eqref{eq:SplitSol}, we then consider solving the 
following variational problem: \textit{Find $(u^0,v^0)\in 
L^2(0,T;H^1_{0,D}(\Omega))\times L^2(0;T;H^1_{0,D}(\Omega))$ such that} 
\begin{equation}
u^0(0)=u_0-u^D(0)\,, \qquad v^0(0)=v_0-v^D(0)
\end{equation}
\textit{and, for all $(\phi,\psi)\in (L^2(0;T;H^1_{0,D}(\Omega)))^2$,}
\begin{align}
\label{eq:cont_weak_formulation_1}
\int_{I} \bigl< \partial_t \displacementC^0, \phi \bigr>_{\Omega}
-
\bigl< \velocityC^0, \phi \bigr>_{\Omega} \, \d t
={}&
0\,, \\
\begin{split}
\label{eq:cont_weak_formulation_2}
\int_{I} \bigl< \partial_t \velocityC^0, \psi \bigr>_{\Omega} +
\langle c^2 \nabla \displacementC^0, \nabla \psi \rangle_{\Omega} \, \d t
={}&
\int_{I} \Big( \bigl< \forceC, \psi \bigr>_{\Omega}
+
\bigl< \partial_n u, \psi \bigr>_{\partial \Omega_N}
\\[1ex]
&\quad -
\bigl< \partial_t \velocityC^D, \psi \bigr>_{\Omega}
-
\langle c^2 \nabla \displacementC^D, \nabla \psi \rangle_{\Omega}\Big)
\, \d t\,.
\end{split}
\end{align}

\begin{remark} i) We note that the correct treatment of inhomogeneous time-dependent 
boundary conditions is an import issue in the application of variational space-time 
methods. The space-time discretization that is derived below (cf.\ \cref{sec:Gal-col_scheme}) and based on the variational problem
\eqref{eq:cont_weak_formulation_1}, \eqref{eq:cont_weak_formulation_2} ensures 
convergence rates of optimal order in space and time, also for time-dependent boundary 
conditions. This is confirmed by the second of the numerical experiments given in \cref{sec:convergence_test_C1}. \\
ii) Our Galerkin--collocation approach is based on solving, along with some collocation 
conditions, the variational equations \eqref{eq:cont_weak_formulation_1}, 
\eqref{eq:cont_weak_formulation_2} in finite dimensional subspaces. In 
particular, the same approximation space will be used for $u^0$ and $v^0$. For this 
reason, the solution space for $v^0$ and the test space in \cref{eq:cont_weak_formulation_2} are chosen slightly stronger than usually; cf.\ 
\cite{Bangerth2010}. 
Choosing $L^2(0;T;L^2(\Omega)))$ instead, would have been sufficient. 
\end{remark}

\section{Galerkin-collocation schemes}
\label{sec:Gal-col_scheme}

In this section we introduce two families of Galerkin--collocation schemes. These 
families combine the concept of collation condition methods applied to the spatially 
discrete counterpart of the equations \eqref{eq:InitialProblem} with the 
finite element discretization of the variational equations 
\eqref{eq:cont_weak_formulation_1}, \eqref{eq:cont_weak_formulation_2}. The collocation 
constrains then allow us to reduce the size of the discrete test spaces for the 
variational conditions compared to a standard Galerkin--Petrov approach; cf.\ 
\cite{Bause2014}. 

First, we need some notation. For the time discretization we decompose the time interval 
$I=(0,T]$ into $N$ subintervals 
$I_n=(t_{n-1},t_n]$, where $n\in \{1,\ldots ,N\}$ and $0=t_0<t_1< \cdots < t_{n-1} < t_n 
= T$ such that $I=\bigcup_{n=1}^N I_n$. We put $\tau = \max_{n=1,\ldots N} \tau_n$ with 
$\tau_n = t_n-t_{n-1}$. Further, the set of time intervals $\mathcal M_\tau := 
\{I_1,\ldots, I_n\}$ is called the time mesh. For a Banach space $B$ and any $k\in \N$, 
we let 
\begin{equation}
\label{Def:Pk}
\mathbb P_k(I_n;B) = \bigg\{w_\tau : I_n \mapsto B \; \Big|\; w_\tau(t) = \sum_{j=0}^k 
W^j t^j 
\,, \; \forall t\in I_n\,, \; W^j \in B\; \forall j \bigg\}\,.
\end{equation}
For an integer $k\in \N$ we introduce the space of globally continuous functions 
in time  
\begin{equation}
\label{Eq:DefXk} 
X_\tau^k (B) := \left\{w_\tau \in C(\overline I;B) \mid w_\tau{}_{|I_n} \in \mathbb 
P_k(I_n;B)\; \forall I_n\in \mathcal M_\tau \right\}\,,
\end{equation}
and for an integer $l\in \N_0$ the space of globally $L^2$-functions in time
\[
\label{Eq:DefYk}
Y_\tau^{l} (B) := \left\{w_\tau \in L^2(I;B) \mid w_\tau{}_{|I_n} \in \mathbb 
P_{l}(I_n;B)\; \forall I_n\in \mathcal M_\tau \right\}\,.
\]

For the space discretization, let $\mathcal T_h=\{K\}$ be a shape-regular mesh of 
$\Omega$  
consisting of quadrilateral or hexahedral elements with mesh size $h>0$. Further, for 
some integer $p\in \N$ let $V_h=V_h^{(p)}$ be the finite element space that is given by 
\begin{equation}
\label{Eq:DefVh}
V_h = V_h^{(p)}=\left\{v_h \in C(\overline \Omega)
\mid v_h{}_{|T} \circ T_K \in \mathbb Q_p \, 
\forall K 
\in \mathcal \Tau_h \right\}\cap H^1_{0,D}(\Omega)\,,
\end{equation}
where $T_K$ is the invertible mapping from the reference cell $\hat K$ to the cell $K$ of 
$\Tau_h$ and $\mathbb Q_p$ is the space of all polynomials of maximum degree $p$ in 
each variable. We let $\mathcal A_h: 
H^1_{0,D}(\Omega) \mapsto V_h$ be the discrete operator that is defined by 
\begin{equation}
\label{Eq:DefAh}
\langle \mathcal A_h w , v_h \rangle = \langle \nabla w, \nabla v_h\rangle  \quad 
\text{for all}\;  v_h\in V_h\,.
\end{equation}
Moreover, $(u_{0,h},v_{0,h})\in V_h^2$ and $(u^D_{\tau,h},v^D_{\tau,h})\in 
(C([0,T];V_h))^2$ define suitable finite element approximations of the initial values 
$(u_0,v_0)$ and the extended boundary values $(u^D,v^D)$ in \cref{eq:SplitSol}. Here, we use interpolation in $V_h$ of the given data.

Now we define our classes of Galerkin--collocation schemes. We follow the 
lines of \cite{Bause2019,Becher2018}.  We restrict ourselves to the schemes studied in 
the numerical experiments presented in Secs.\ 
\ref{sec:convergence_test_C1}, \ref{sec:Challenging_Example} and \ref{sec:NumTestC2}. The 
definition of classes of  Galerkin--collocation 
schemes 
with even higher regularity in time is straightforward, but not done here.

\begin{definition}[$\boldsymbol{C^l}$--regular in time Galerkin-collocation schemes GCC$\boldsymbol{ 
		{}^l(k)}$] 
	\label{Def:GCC}
	Let $l\in \{1,2\}$ be fixed and $k\geq 2l+1$. For $n=1,\ldots, N$ and 
	given $(u_{\tau,h}{}_|{}_{I_{n-1}}(t_{n-1})$, $v_{\tau,h}{}_|{}_{I_{n-1}}(t_{n-1}))\in 
	V_h^2$  for $n>1$ and $u_{\tau,h}{}_|{}_{I_0}(t_0)=u_{0,h}$, 
	$v_{\tau,h}{}_|{}_{I_0}(t_0)=v_{0,h}$ for $n=1$,  find $(u^0_{\tau,h}{}_|{}_{I_n}, 
	v^0_{\tau,h}{}_|{}_{I_n})  \in 
	(\mathbb P_k (I_n;V_h))^2$ such that, for $s_0\in \N_0$, $s_1\in \N$ with $s_0,s_1\leq l$, 
	\begin{align}
	%
	& \partial_t^{s_0} w^0_{\tau,h}{}_|{}_{I_n}(t_{n-1})
	= \partial_t^{s_0} w^0_{\tau,h}{}_|{}_{I_{n-1}}(t_{n-1})
	\,, \quad \text{for} \; w^0_{\tau,h} \in \left\{u^0_{\tau,h},v^0_{\tau,h}\right\}\,,
	\label{Eq:SemiDisLocalcGPC_6}
	\\[1ex]
	&\partial_t^{s_1} u^0_{\tau,h}{}_|{}_{I_n}(t_{n})
	-\partial_t^{s_1-1}v^0_{\tau,h}{}_|{}_{I_n}(t_{n})
	= 0 \,,
	\label{Eq:SemiDisLocalcGPC_7}
	\\[1ex]
	%
	&
	\begin{aligned}
	\partial_t^{s_1} v^0_{\tau,h}{}_|{}_{I_n}(t_{n}) 
	+ \mathcal A_h \partial_t^{s_1-1} & u^0_{\tau,h}{}_|{}_{I_n}(t_{n})
	= \partial_t^{s-1} f(t_{n})\\
	& - \partial_t^{s_1} 
	v^D_{\tau,h}{}_|{}_{I_n}(t_{n}) - \mathcal A_h \partial_t^{s_1-1} u^D_{\tau,h}{}_|{}_{I_n}(t_{n}) 
	\,,
	\label{Eq:SemiDisLocalcGPC_8}
	\end{aligned}
	\end{align}
	and, for all $(\varphi_{\tau,h}, \psi_{\tau,h}) \in 
	(\mathbb P_{0} (I_n;V_h))^2$, 
	\begin{align}
	\label{Eq:SemiDisLocalcGPC_9}
	&\int_{I_n} \Big(\langle \partial_t u^0_{\tau,h} , \varphi_{\tau,h} \rangle_{\Omega} 
	- \langle v^0_{\tau,h} , \varphi_{\tau,h} \rangle_{\Omega} \Big) \, \d t
	= 
	0\,,
	\\
	& \begin{alignedat}{1}
	\int_{I_n} \Big(\langle \partial_t v^0_{\tau,h} , \psi_{\tau,h} \rangle_{\Omega} 
	& + \langle \mathcal A u^0_{\tau,h} , \psi_{\tau,h} \rangle_{\Omega} \Big) \, \d t
	= 
	\int_{I_n}
	\langle f,\psi_{\tau,h}\rangle_{\Omega} \, \d t
	\\
	& -
	\int_{I_n} \Big(\langle \partial_t v^D_{\tau,h} , \psi_{\tau,h} \rangle_{\Omega} 
	+ \langle \mathcal A_h u^D_{\tau,h} , \psi_{\tau,h} \rangle_{\Omega} \Big) \, \d t\,.
	\label{Eq:SemiDisLocalcGPC_10}
	\end{alignedat}
	\end{align}
	
\end{definition}

\begin{remark}
\begin{itemize}
\item In \cref{Eq:SemiDisLocalcGPC_6}, the discrete 
initial values $(\partial_t u_{\tau,h}(0),\partial_t v_{\tau,h}(0))$ arise for $s_0=1$. For 
$\partial_t u_{\tau,h}(0)$ we use a suitable finite element approximation $v_{0,h}\in 
V_h$ (here, an interpolation) of $v_0\in V$. For $\partial_t v_{\tau,h}(0)$ 
we evaluate the wave equation in the initial time point and use a suitable finite element 
approximation (here, an interpolation) of $\partial_t^2 u(0)= c^2 \Delta u(0) + f(0)$. For $s_0=2$ in \cref{Eq:SemiDisLocalcGPC_6}, the initial value $\partial_t^2 
v_{\tau,h}(0)$ is computed as a suitable finite element approximation (here, an 
interpolation) of $\partial_t^3 u(0)= c^2 \Delta \partial_t u(0) + \partial_t f(0)$. 
Mathematically, this approach requires that the partial equation and its time derivative 
are satisfied up to the initial time point and, thereby, sufficient regularity of the 
continuous solution. Without such regularity assumptions, the application of higher order 
discretization schemes cannot be justified rigorously. Nevertheless, in practice such 
methods often show a superiority over lower-order ones, even for solutions without the 
expected high regularity (cf.\ \cref{sec:Challenging_Example}).        
 
\item From \cref{Eq:SemiDisLocalcGPC_6}, $(u_{\tau,h},v_{\tau,h})\in 
(C^l(\overline I;V_h))^2$, for fixed $l \in \{1,2\}$, is easily concluded. 
\end{itemize}
\end{remark}

An optimal order error analysis for the GCC$^1(k)$ family of schemes of Def.\ 
\ref{Def:GCC} is provided in \cite{Bause2019}. The following theorem is proved.
\begin{theorem}[Error estimates for ${\boldsymbol{(u_{\tau,h}$, $v_{\tau,h})}}$ of 
GCC$^1(\boldsymbol k)$]
	\label{th:error}
	Let $l=1$ and $k\geq 3$. For the error $(e^{\; u}, e^{\; v})=(u-u_{\tau,h}, v-v_{\tau,h})$ of the fully discrete scheme GCC$^l(k)$ of Def.~\ref{Def:GCC} there holds that 
	\begin{equation}
	\begin{aligned}
	\label{eq:error_1a}
	\| e^{\; u}(t)\| + \| e^{\; v}(t)\|
	&\lesssim
	\tau^{k+1}+h^{p+1} \, , \;\; t\in \overline I\,,
	\\
	\|\nabla e^{\; u}(t)\| 
	&\lesssim
	\tau^{k+1}+h^{p}  \, , \;\; t\in \overline I\,,
	\end{aligned}
	\end{equation}
	as well as
	\begin{equation}
	\begin{aligned}
	\label{eq:error_2a}
	\| e^{\; u}(t)\|_{L^2(I;H)} + \| e^{\; v}(t)\|_{L^2(I;H)}
	&\lesssim   \tau^{k+1}+h^{p+1} \, ,
	\\
	\|\nabla e^{\; u}(t)\|_{L^2(I;H)}
	&\lesssim
	\tau^{k+1}+h^{p} \, .
	\end{aligned}
	\end{equation}
\end{theorem}

Error estimates for the GCC$^2(k)$ family remain as a work for the future. In \cref{sec:NumTestC2}, the convergence of GCC$^2(5)$ is demonstrated numerically. 
Further, we note that a computationally cheap post-processing of improved regularity and  
accuracy for continuous Galerkin--Petrov methods is presented and studied in 
\cite{Bause2018}.


In the next sections we study the schemes GCC$^1(3)$ and GCC$^2(5)$ of 
Def.~\ref{Def:GCC} in detail. Their algebraic forms are derived and the algebraic linear 
solver are presented. Finally, the results of our numerical experiments with the proposed 
methods are presented. Here, we restrict ourselves to the lowest-order cases with $k=3$ 
for $l=1$ and $k=5$ for $l=2$ of Def.~\ref{Def:GCC}. This is sufficient to demonstrate the 
potential of the Galerkin--collocation approach and its superiority over the standard 
continuous Galerkin approach in space and time \cite{KM04,Bause2018}. An implementation of 
GCC$^l(k)$ for higher values of $k$ along with efficient algebraic solvers is currenty 
still missing.

\section{\texorpdfstring{Galerkin--collocation GCC$\boldsymbol{{}^1(3)}$}
	{Galerkin--collocation GCC1(3)}}
\label{sec:c1_solution}

Here, we derive the algebraic system of the GCC$^1(3)$ approach and discuss our algebraic solver for the arising block system. For brevity, the derivation is done for $k=3$ only. The generalization to larger values of $k$ is straightforward; cf.\ \cite{Bause2019}.

\subsection{Fully discrete system}
\label{sec:Deriving_C1_System}

To derive the discrete counterparts of the variational conditions \eqref{Eq:SemiDisLocalcGPC_9}, \eqref{Eq:SemiDisLocalcGPC_10} and the collocation constrains  \eqref{Eq:SemiDisLocalcGPC_6} to \eqref{Eq:SemiDisLocalcGPC_8}, we let $\{\phi_j\}_{j=1}^J \subset V_h$, denote a (global) nodal Lagrangian basis of $V_h$. The mass matrix $\mat M$ and the stiffness matrix $\mat A$ are defind by 
\begin{align}
\label{eq:DefM_A}
	\mat{M} &\coloneqq
	 \left(\bigl<
	 \phi_i, \phi_j
	 \bigr>_{\Omega}\right)_{i,j=1}^J\,,
	 &
 	\mat{A} &\coloneqq
	 \left(\bigl<
	 \nabla \phi_i, \nabla \phi_j
	 \bigr>_{\Omega}\right)_{i,j=1}^J\,,
\end{align}
On the reference time interval $\hat I = [0,1]$ we define a Hermite-type basis $\{\hat 
\xi_l\}_{l=0}^3 \subset \mathbb P_3 (\hat I;\R)$ of $\mathbb P_3 (\hat I;\R)$ by the 
conditions 
\begin{equation}
\label{eq:CondXi}
\begin{aligned}
\hat{\xi_{0}}(0) &= 1\,,
& \hspace{1em}
\hat{\xi_{0}}(1) &= 0\,,
& \hspace{1em}
\partial_t\hat{\xi_{0}}(0) &= 0\,,
& \hspace{1em}
\partial_t\hat{\xi_{0}}(1) &= 0\,,
\\
\hat{\xi_{1}}(0) &= 0\,,
&
\hat{\xi_{1}}(1) &= 0\,,
&
\partial_t\hat{\xi_{1}}(0) &= 1\,,
&
\partial_t\hat{\xi_{1}}(1) &= 0\,,
\\
\hat{\xi_{2}}(0) &= 0\,,
&
\hat{\xi_{2}}(1) &= 1\,,
&
\partial_t\hat{\xi_{2}}(0) &= 0\,,
&
\partial_t\hat{\xi_{2}}(1) &= 0\,,
\\
\hat{\xi_{3}}(0) &= 0\,,
&
\hat{\xi_{3}}(1) &= 0\,,
&
\partial_t\hat{\xi_{3}}(0) &= 0\,,
&
\partial_t\hat{\xi_{3}}(1) &= 1\,.
\end{aligned}
\end{equation}
\
These conditions then define the basis of $\mathbb P_3(\hat I;\R)$ by  
\begin{align}
\label{eq:basisxi}
\hat{\xi_{0}} &= 1 - 3 t^2 + 2 t^3, &
\hat{\xi_{1}} &= t - 2 t^2 + t^3, &
\hat{\xi_{2}} &= 3 t^2 - 2 t^3, &
\hat{\xi_{3}} &= -t^2 + t^3.
\end{align}
By means of the affine transformation $T_n(\hat t):= t_{n-1} + \tau_n\cdot \hat t$, 
with $\hat t  \in \hat I$, from the reference interval $\hat I$ to $I_n$ such that 
$t_{n-1}=T_n(0)$ and $t_{n}=T_n(1)$, the basis $\{\xi_l\}_{l=0}^3 \subset \mathbb P_3 
(I_n;\R)$ is given by $\xi_l = \hat \xi_l \circ T_n^{-1}$ for $l=0,\ldots, 3$. In terms of 
basis functions, $w_{\tau,h}\in \mathbb P_3 (I_n;V_h)$ is thus represented by 
\begin{align}
\label{eq:full_discrete_ansatz}
w_{\tau,h}(\mat x, t) = \sum_{l=0}^{3} \sum_{j=1}^{J}w_{n,l,j} \phi_j 
(\mat{x})\xi_{l}(t)\,, \quad (\mat x,t)\in \Omega \times \overline{I_n}\,.
\end{align}
For $\zeta_0 \equiv 1$ on $\overline{I_n}$, a test basis of $\mathbb P_0 (I_n;V_h)$ is then given by 
\begin{equation}
 \label{eq:testbasis}
\mathcal B = \{\phi_1 \zeta_0, \ldots, \phi_J \zeta_0\}\,.
\end{equation}
To evaluate the time integrals on the right-hand side of \cref{Eq:SemiDisLocalcGPC_10} 
we still use the Hermite-type interpolation operator $I_{\tau|I_n}$, on $I_n$, defined 
by 
\begin{equation}
\label{eq:hermite_interpolation}
\begin{alignedat}{2}
I_{\tau|I_n}g(t)
&\coloneqq
 \sum_{s=0}^{l} \tau_n^s \hat \xi_s  (0)
\underbrace{\partial_t^s g_|{}_{I_n}(t_{n-1})}_{=:g_s}
&&+ \sum_{s=0}^{l} \tau_n^s \hat \xi_{s+l+1} (1)
\underbrace{\partial_t^s g|_{I_n}(t_{n})}_{=:g_{s+l+1}}\,.
\end{alignedat}
\end{equation}
Here, the values $\partial_t^s g_|{}_{I_n}(t_{n-1})$ and 
$\partial_t^s g_|{}_{I_n}(t_{n})$ in \eqref{eq:hermite_interpolation} denote the corresponding 
one-sided limits of values $\partial_t^s g(t)$ from the interior of $I_n$.

Now, we can put the equations of the proposed GCC$^1(3)$ approach in their algebraic 
forms. In the variational equations \eqref{Eq:SemiDisLocalcGPC_9} and 
\eqref{Eq:SemiDisLocalcGPC_10}, we use the representation 
\eqref{eq:full_discrete_ansatz} for each component of $(u_{\tau,h}, v_{\tau,h}) 
\in (\mathbb P_3(I_n;V_h))^2$, choose the test functions \eqref{eq:testbasis} 
and interpolate the right-hand side of \eqref{Eq:SemiDisLocalcGPC_10} by applying 
\eqref{eq:hermite_interpolation}. All of the arising time integrals are evaluated 
analytically. Then, we can recover the variational conditions 
\eqref{Eq:SemiDisLocalcGPC_9} and \eqref{Eq:SemiDisLocalcGPC_10} on the subinterval $I_n$ 
in their algebraic forms
\begin{align}
\label{eq:time_scheme_a1}
&\mat{M} \left(
-\mat{u}_{n,0}^{0} + \mat{u}_{n,2}^{0}
\right)
-
\tau_n \mat{M} \left(
\frac{1}{2} \mat{v}_{n,0}^{0} + \frac{1}{12} \mat{v}_{n,1}^{0} + \frac{1}{2} 
\mat{v}_{n,2}^{0} 
- \frac{1}{12} \mat{v}_{n,3}^{0}
\right)
=
\mat 0\,,
\\
\begin{split}
\label{eq:time_scheme_a2}
&\mat{M} \left(
-\mat{v}_{n,0}^{0} + \mat{v}_{n,2}^{0}
\right)
+
\tau_n \mat{A}
\left(
\frac{1}{2} \mat{u}_{n,0}^{0} + \frac{1}{12} \mat{u}_{n,1}^{0} + \frac{1}{2} 
\mat{u}_{n,2}^{0} 
- \frac{1}{12} \mat{u}_{n,3}^{0}
\right)
=
\\
&
\hphantom{
	\mat{M} \left(
	-\mat{v}_{n,0}^{D}\right)
}
\tau_n \mat{M}
\left( \frac{1}{2} \mat{f}_{n,0}^{} + \frac{1}{12} \mat{f}_{n,1}^{} + \frac{1}{2} 
\mat{f}_{n,2}^{} - \frac{1}{12} \mat{f}_{n,3} \right)
-
\mat{M} \left(
-\mat{v}_{n,0}^{D} + \mat{v}_{n,2}^{D}
\right)
\\
&
\hphantom{
	\mat{M} \left(
	-\mat{v}_{0}^{D} + \mat{v}_{2}^{K,I} + \mat{v}_{2} +\mat{v}_{2}
	\right)
}
-
\tau_n \mat{A}
\left(
\frac{1}{2} \mat{u}_{n,0}^{D} + \frac{1}{12} \mat{u}_{n,1}^{D} + \frac{1}{2} 
\mat{u}_{2}^{D} 
- \frac{1}{12} \mat{u}_{n,3}^{D}
\right).
\end{split}
\end{align}
This gives us the first two equations for the set of eight unknown solution vectors 
$\mathcal L = \{\mat u_{n,0}^0,\ldots, \mat u_{n,3}^0,\mat v_{n,0}^0,\ldots, \mat 
v_{n,3}^0\}$ on each subinterval $I_n$, where each of these vectors is 
defined by means of \eqref{eq:full_discrete_ansatz} through $\mat w = 
(w_1,\ldots,w_J)^\top$ for $\mat w \in \mathcal L$.

Next, we study the algebraic forms of the collocations conditions 
\eqref{Eq:SemiDisLocalcGPC_6}  to \eqref{Eq:SemiDisLocalcGPC_8}. By means 
of the definition \eqref{eq:CondXi} of the basis of $\mathbb P(I_n;\R)$, the  
constraints \eqref{Eq:SemiDisLocalcGPC_6} read as 
\begin{equation}
\label{eq:constraint_left_a}
\mat{u}_{n,0}^{0} = \mat{u}_{n-1,2}^{0}\,, \;\; \mat{u}_{n,1}^{0} = 
\mat{u}_{n-1,3}^{0}\,, \quad \mat{v}_{n,0}^{0} = \mat{v}_{n-1,2}^{0}\,, \;\; 
\mat{v}_{n,1}^{0} = \mat{v}_{n-1,3}^{0}\,.
\end{equation}
By means of \eqref{eq:CondXi} along with \eqref{eq:DefM_A}, the conditions 
\eqref{Eq:SemiDisLocalcGPC_7} and \eqref{Eq:SemiDisLocalcGPC_8} can be recovered as 
\begin{align}
\label{eq:C1_Col2}
\mat{M}
\frac{1}{\tau_n} \mat{u}_{n,3}^{0}
-
\mat{M}
\mat{v}_{n,2}^{0}
& = \mat 0\,,
\\
\label{eq:C1_Col2_2}
\mat{M} \frac{1}{\tau_n} \mat{v}_{n,3}^{0}
+
\mat{A} \mat{u}_{n,2}^{0}
& =
\mat{M}
\mat{f}_{n,2}
-
\mat{M} \frac{1}{\tau_n} \mat{v}_{n,3}^{D}
-
\mat{A} \mat{u}_{n,2}^{D}\,.
\end{align}

Putting relations \eqref{eq:constraint_left_a} into the identities 
\eqref{eq:time_scheme_a1} and \eqref{eq:time_scheme_a2} and combining the resulting 
equations with \eqref{eq:C1_Col2} and \eqref{eq:C1_Col2_2} yields for the subinterval 
$I_n$ the linear block system 
\begin{equation}
\label{eq:Sx=b}
\mat{S}\mat{x} = \mat{b}
\end{equation}
for the vector of unknowns 
\begin{equation}
\label{def:vec_x}
\mat{x} = 
\left(\left(\mat{v}_{n,2}^{0}\right)^{\top}\,, \left(\mat{v}_{n,3}^{0}\right)^\top\,, 
\left(\mat{u}_{n,2}^{0}\right)^\top\,, \left(\mat{u}_{n,3}^{0}\right)^\top\right)^\top
\end{equation}
and the system $\mat S$ and right-hand side $\mat b$ given by 
\begin{align}
\label{eq:system_matrix_c1}
\mat{S} &=
\begin{pmatrix}
\mat{M} & \mat{0} & \mat{0} & \frac{1}{\tau_n}\mat{M} \\[1ex]
\mat{0} & \frac{1}{\tau_n}\mat{M} & \mat{A} & \mat{0} \\[1ex]
-\frac{\tau_n}{2}\mat{M} & \frac{\tau_n}{12}\mat{M} & \mat{M} & \mat{0} \\[1ex]
\mat{M} & \mat{0} & \frac{\tau_n}{2}\mat{A} & -\frac{\tau_n}{12}\mat{A}
\end{pmatrix},
&
\mat{b} &= 
\begin{pmatrix}
\mat{0}
\\[1ex]
\mat{M}
\left(\mat{f}_{n,2}
- \frac{1}{\tau_n}  \mat{v}^D_{n,3}
\right)
-
\mat{A} \mat u^D_{n,2}
\\[1ex]
\mat{M}
\left(
\mat{u}^0_{n,0}
+ \frac{\tau_n}{2}\mat{v}^0_{n,0}
+ \frac{\tau_n}{12}\mat{v}^0_{n,1}
\right)
\\[1ex]
\mat{b}_{n,4}
\end{pmatrix},
\end{align}	
with $\mat{b}_{n,4} = \mat{M} \bigl(
\mat{v}^0_{n,0}
+ \mat{v}^D_{n,0}
- \mat{v}^D_{n,2}
+ \frac{\tau_n}{2} (\mat{f}_{n,0} + \mat{f}_{n,2})
+ \frac{\tau_n}{12} (\mat{f}_{n,1} - \mat{f}_{n,3})
\bigr)
- \mat{A}
\bigl(
\frac{\tau_n}{2} (\mat{u}^0_{n,0} + \mat{u}^D_{n,0} + \mat{u}^D_{n,2})
+
\frac{\tau_n}{12} (\mat{u}^0_{n,1} + \mat{u}^D_{n,1} - \mat{u}^D_{n,3})
\bigr)$. By means of the collocation constraints \eqref{eq:constraint_left_a}, the number of unknown coefficient vectors for the discrete 
solution $(u_{\tau,h}{}_{|I_n},v_{\tau,h}{}_{|I_n})\in (\mathbb P_3(I_n;V_h))^2$ is 
thus effectively reduced from eight to four vectors, assembled now in $\mat x$ by \eqref{def:vec_x}. 

We note that the first two rows of \cref{eq:system_matrix_c1} represent the 
collocation conditions \eqref{eq:C1_Col2} and \eqref{eq:C1_Col2_2}. They have a sparser 
structure then the last two rows representing the variational conditions which can be 
advantageous or exploited for the construction of efficient iterative solvers for 
\eqref{eq:Sx=b}. Compared with a pure variational approach (cf.\ 
\cite{Hussain2011,Hussain2011_2,KM04}), more degrees or freedom are obtained directly 
by computaionally cheap vector identities (cf.\ \eqref{eq:time_scheme_a1}) in GCC$^l(k)$ such that they 
can be eliminated from the overall linear system and, thereby, used to reduce the systems size.

\subsection{Solver technology}
\label{sec:solver_technology}

In the sequel, we present two different iterative approaches for solving the linear 
system \eqref{eq:Sx=b} with the non-symmetric matrix $\mat S$. In \cref{sec:Challenging_Example}, a runtime comparison between the two concepts is 
provided. As basic toolbox we use the deal.II finite element 
library \cite{DealIIReference} along with the Trilinos library \cite{TrilinosReference} for parallel 
computations. 

\subsubsection{1.\ Approach: Condensing the linear system}
\label{sec:method_1}
The first method for solving \eqref{eq:Sx=b} is based on the concepts developed in 
\cite{Koecher15}. The key idea is to use Gaussian block elimination within the system 
matrix $\mat{S}$ and end up with a linear system with matrix $\mat{S}_r$ of reduced size 
for one of the subvectors in $\mat x$ in \eqref{def:vec_x} only, and to compute the 
remaining subvectors of \eqref{def:vec_x} by computationally cheap post-processing 
steps afterwards. The reduced system matrix $\mat{S}_r$ should have sufficient potenial that  
an efficient preconditioner for the iterative solution of the reduced 
system can be constructed. Of course, the Gauss elimination on the block level can be done in different 
ways. The goal of our approach is to avoid the inversion of the stiffness matrix 
$\mat A$ in \eqref{eq:system_matrix_c1} in the computation of the condensed system 
matrix $\mat{S}_r$ such that a matrix-vector multiplication with $\mat{S}_r$ just involves 
calculating $\mat{M}^{-1}$. At least for discontinuous Galerkin methods in space, where $\mat{M}$ is block diagonal, this is computationally 
cheap; cf.\ \cite{Bause2014,Koecher15}. We note that a continuous Galerkin approach in 
space is used here only in order to simplify the notation and since the discretization in time 
by the combined Galerkin--collocation approach is in the scope of interest. 

Here, we choose the subvector $\mat u_{n,2}^0$ of $\mat x$ in \eqref{def:vec_x} as the 
essential unknown, i.e.\ as the unknown solution vector of the condensed system with  
matrix $\mat{S}_r$. By block Gaussian elemination we then end up with solving the linear 
system, 
\begin{equation}
\label{eq:recuced_c1_system}
\left(
\mat{M} +  \frac{\tau_n^2}{12} \mat{A} + \frac{\tau_n^4}{144} \mat{A}\mat{M}^{-^1}\mat{A}
\right)
\mat{u}_{n,2}^{0}
= \mat{b}_{n,r}
\end{equation}
with right-hand side vector
\begin{multline}
\label{def:br}
\mat{b}_{n,r} =
\mat{M} \Biggl(
\frac{1}{2} \mat{f}_{n,0} + \frac{1}{12} \mat{f}_{n,1} + \frac{1}{3} \mat{f}_{n,2} - \frac{1}{12} 
\mat{f}_{n,3}
\Biggr)
+ \mat{M} \Biggl(
2 \mat{v}_{n,0}^0 + \frac{1}{6} \mat{v}_{n,1}^0 + \frac{2}{\tau_n} \mat{u}_{n,0}^0
\Biggr)
\\
- \mat{A} \Biggl(
\frac{2}{3} \tau_n \mat{u}_{n,0}^0 \! + \! \frac{1}{12} \tau_n \mat{u}_{n,1}^0 \!+ \!\frac{1}{12} \tau_n ^2 
\mat{v}_{n,2}^0 \!+ \!\frac{1}{72} \tau_n^2 \mat{v}_{n,1}^0
\Biggr)
\!+ \! \mat{M} \Biggl(
2 \mat{v}_{n,0}^D \! + \!\frac{1}{6} \mat{v}_{n,1}^D \!+ \!\frac{2}{\tau_n} \mat{u}_{n,0}^D \!- \! \frac{2}{\tau_n} 
\mat{u}_{n,2}^D
\Biggr)
\\
+\mat{A} \Biggl(
\frac{1}{72} \tau_n^3 \mat{f}_{n,2} -  \mat{M}^{-1} \frac{1}{72} \tau_n^3 \mat{u}_{n,2}^D
\Biggr)
\! + \!
\mat{A} \Biggl(
- \frac{2}{3} \tau_n \mat{u}_{n,0}^D \!- \! \frac{1}{12} \tau_n \mat{u}_{n,1}^D \! -\!  \frac{\tau_n}{6} 
\mat{u}_{n,2}^D \! - \! \frac{\tau_n^2}{12} \mat{v}_{n,0}^D \! -\!  \frac{\tau_n^2}{72} \mat{v}_{n,1}^D
\Biggr)\,.
\end{multline}
The product of $\mat A \mat{M}^{-1}\mat{A}$ in \eqref{eq:recuced_c1_system} mimics the 
discretization of a fourth order operator due to the appearance of the product of $\mat A$ 
with its ''weighted'' form $\mat{M}^{-1}\mat{A}$. Thereby, the conditioning number of the condensed system is strongly 
increased (cf.\ \cite{Koecher15}) which is the main drawback in this concept of condensing the 
overall system \eqref{eq:Sx=b} to \eqref{eq:recuced_c1_system} for the essential unknown 
$\mat u_{n,2}^0$. On the other hand, since $\mat{M}$ and $\mat{A}$ are symmetric and, 
thus, $\mat{A}\mat{M}^{-1}\mat{A} = (\mat{A}\mat{M}^{-1}\mat{A})^\top$, the condensed 
matrix $\mat{S}_r$ is symmetric such that the preconditioned conjugate gradient method 
can be applied. Solving systems of type \eqref{eq:recuced_c1_system} is carefully studied 
in \cite{Bause2014,Koecher15} and the references given therein. 

\begin{figure}[!htb]
\centering
\includegraphics[width=1\textwidth,height=0.8\textheight,keepaspectratio]
{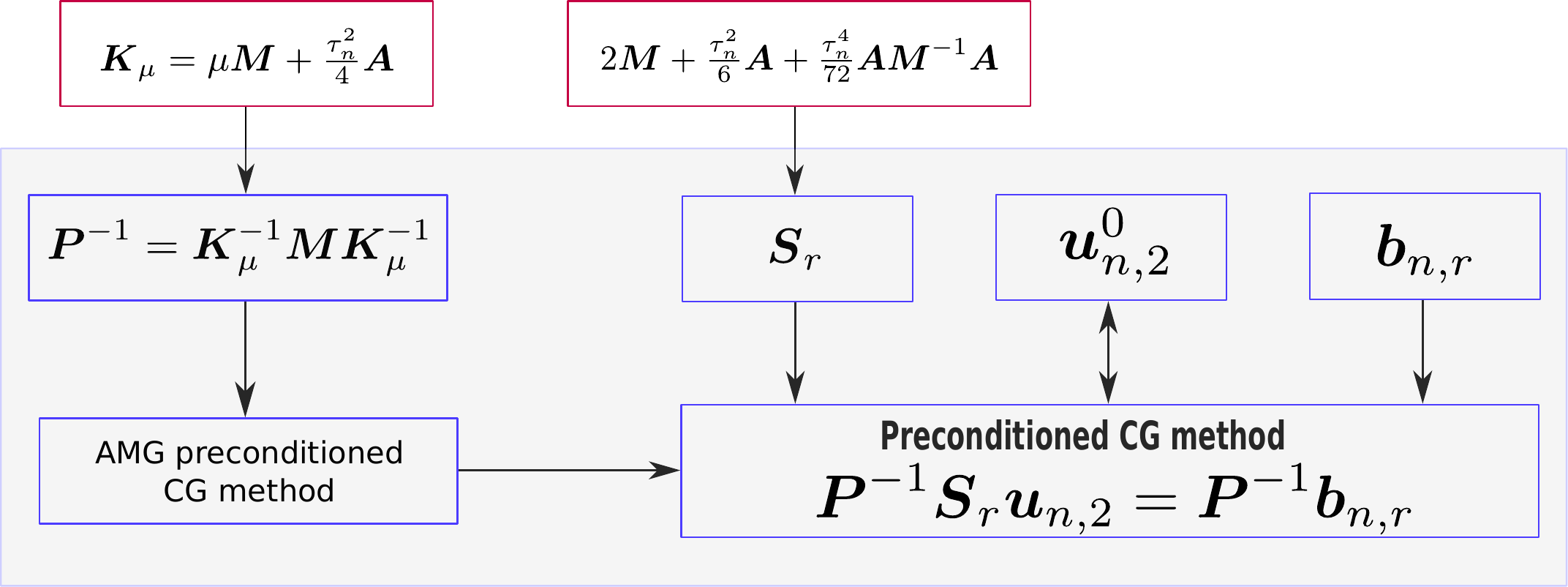}
	\caption{Preconditioning and solver for the condensed system \eqref{eq:recuced_c1_system} of GCC$^1(3)$.}
	\label{fig:Solver_condensed}
\end{figure}

We solve \eqref{eq:recuced_c1_system} by the conjugate gradient method. The left
preconditioning operator
\begin{align}
\mat{P} &= \mat{K}_{\mu} \mat{M}^{-1} \mat{K}_{\mu}
= \left(\mu \mat{M} + \frac{\tau_n^2}{4} \mat{A}\right) \mat{M}^{-1}
\left(\mu \mat{M} + \frac{\tau_n^2}{4} \mat{A}\right)\,,
\end{align}
with positive $\mu \in \R$ chosen such that the spectral norm of 
$\mat{P}^{-1}\mat{S}_r$ is 
minimised, is applied. For details 
of the choice of the parameter $\mu$, we refer to \cite{Bause2014,Koecher15}. Here, we 
use $\mu = \sqrt{11/2}$. In order to apply the preconditioning operator $\mat P$ in the 
conjugate gradient iterations, without assembling $\mat P$ explicitly,  i.e.\ to solve the auxiliary system with matrix $\mat P$,  
we have to solve linear systems for the mass matrix $\mat{M}$ and the stiffness 
matrix $\mat{A}$. For this, we use embedded conjugate gradient iterations combined with 
an algebraic multigrid preconditioner of the Trilinos library \cite{TrilinosReference}. 
The overvall algorithm for solving \eqref{eq:recuced_c1_system} is sketched in \cref{fig:Solver_condensed}. The advantage of this approach is that we just have to store 
$\mat{M}$ and $\mat{A}$ as sparse matrices in the computer memory. We never have to 
assemble the full matrix $ \mat S$ from \eqref{eq:system_matrix_c1}, nor do we have to store 
the reduced matrix $\mat{S}_r$ from \eqref{eq:recuced_c1_system}. Finally, the 
remaining unknown subvectors $\mat{v}_{n,2}^{0}, \mat{v}_{n,3}^{0}$ and 
$\mat{u}_{n,3}^{0}$ in \eqref{eq:Sx=b} are successively computed in post-processing steps. 

\subsubsection{2.\ Approach: Solving the non-symmetric system}
\label{sec:method_2}
The second approach used to solve \eqref{eq:Sx=b} relies on assembling the system 
matrix $\mat{S}$ of \eqref{eq:system_matrix_c1} as a sparse matrix and solving the 
resulting non-symmetric system. For smaller dimensions of $\mat S$ a parallel direct 
solver \cite{Demmel2003} is used. For constant time step sizes $\tau_n$ the matrix $\mat S$ needs to be 
factorized once only, which results in excellent performance properties for large 
sequences of time steps. For high-dimensional problems with interest in practice, 
we use the Generalized Minimal Residual (GMRES) method, an iterative Krylov subspace 
method, to solve \eqref{eq:Sx=b}. The drawback of this approach then comes through the 
necessity to provide an efficient preconditioner, i.e.\  an approximation to the inverse 
of $\mat{S}$, for the complex block matrix $\mat S$ of \eqref{eq:system_matrix_c1}. 
Here, we use the algebraic multigrid method as preconditioning technique. We use the 
MueLue preconditioner \cite{MueLueReference}, which is part of the Trilinos project, with 
non-symmetric smoothed aggregation. We use an usual V-cycle algorithm along with a 
symmetric successive over-relaxation (SSOR) smoother with 
a damping factor of $1.33$. The design of efficient 
algebraic solvers for block systems like \eqref{eq:Sx=b}, and for higher order 
variational time discretizations in general, is still an active field of research. We 
expect further improvement in the future.  

\subsection{Numerical convergence tests}
\label{sec:convergence_test_C1}

In this section we present a numerical convergence test for the proposed GCC$^1(3)$ approach of 
Def.\ \ref{Def:GCC} and \cref{sec:Deriving_C1_System}, respectively. For the 
solution $\{u, v\}$ of \cref{eq:InitialProblem} and the fully discrete  
approximation GCC$^1(3)$ of Def.\ \ref{Def:GCC} we let 
\begin{align}
e^{u} &\coloneqq u(\mat{x},t) - u_{\tau,h}(\mat{x},t),
&
e^{v} &\coloneqq v(\mat{x},t) - v_{\tau,h}(\mat{x},t).
\end{align}
We study the error $(e^u,e^v)$ with respect to the norms 
\begin{equation}
\begin{aligned}
\| e^w \|_{L^\infty(L^2)} &\coloneqq  \max_{t \in I} \left( \int_{\Omega} \| e^w \|^2 \d 
x \right)^{\frac{1}{2}} \,,
\quad
\| e^w \|_{L^2(L^2)} \coloneqq  \biggl( \int_{I} \int_{\Omega} \| e^w \|^2 \d x  \, \d t 
\biggr)^{\frac{1}{2}}\,,
\end{aligned}
\end{equation}
where $w \in (u, v)$, and in the energy quantities
\begin{equation}
\begin{aligned}
||| E \, |||_{L^\infty} &\coloneqq \max_{t \in I} \biggl(\| \nabla e^{u} \|^2 + \| e^{v} 
\|^2 
\biggr)^{\frac{1}{2}} \,,
\quad
||| E\, |||_{L^2} \coloneqq \Biggl( \int_I
\int_{\Omega} \| \nabla e^{u} \|^2 + \| e^{v} \|^2 \d \mat{x} \, \d t
\Biggr)^{\frac{1}{2}} \, .
\end{aligned}
\end{equation}
All $L^{\infty}$-norms in time are computed on the discrete time grid
\begin{equation}
\label{eq:discrete_time_grid}
I = \big\{
t_n^d | t_n^d = t_{n-1} + d \cdot k_n \cdot \tau_n,
\quad
k_n = 0.001, d = 0, \ldots, 999, n = 1, \ldots, N
\big\}.
\end{equation}

For our first convergence test we prescribe the solution 
\begin{equation}
\label{eq:conv_test_1}
u_1(\mat{x},t) = \sin(4 \pi t) \cdot x_1 \cdot (x_1 - 1) \cdot x_2 \cdot (x_2 - 1).
\end{equation}
on $\Omega \times I = (0,1)^2 \times [0,1]$. We let $c = 1$, use a constant mesh size 
$h_0 = 0.25$ and start with the time step size $\tau_0 = 0.1$. We compute the errors on a 
sequence of successively refined time meshes by halving the step sizes in each 
refinement step. We choose a bicubic discretization of the space 
variables in $V_h^3$ (cf.\ \eqref{Eq:DefVh})  such that the spatial part of the solution 
is resolved exactly by its numerical approximation. Table 
\ref{tab:conv_1} summarizes the computed errors and experimental orders of convergence. 
The expected convergence rates of Thm.~\ref{th:error} are nicely confirmed.
\begin{table}[!htb]
	\centering
	\small
	\begin{tabular}{c@{\,\,\,\,}c  c@{\,}c  c@{\,}c  c@{\,}c}
		\toprule
		{$\tau$} & {$h$} &
		{ $\| e^{u}  \|_{L^\infty(L^2)} $ } & EOC &
		{ $\| e^{v}  \|_{L^\infty(L^2)} $ } & EOC &
		{ $||| E\,  |||_{L^\infty} $ } & EOC \\
		\midrule
		$\tau_0/2^0$ & $h_0$ & 2.318e-04 & {--} & 1.543e-03 & {--} & 1.574e-03 & 
{--} \\
		$\tau_0/2^1$ & $h_0$ & 1.541e-05 & 3.91 & 9.694e-05 & 3.99 & 1.004e-04 & 
3.97 \\
		$\tau_0/2^2$ & $h_0$ & 9.825e-07 & 3.97 & 6.260e-06 & 3.95 & 6.478e-06 & 
3.95 \\
		$\tau_0/2^3$ & $h_0$ & 6.185e-08 & 3.99 & 3.946e-07 & 3.99 & 4.082e-07 & 
3.99 \\
		$\tau_0/2^4$ & $h_0$ & 3.873e-09 & 4.00 & 2.472e-08 & 4.00 & 2.557e-08 & 
4.00 \\
		$\tau_0/2^5$ & $h_0$ & 2.422e-10 & 4.00 & 1.548e-09 & 4.00 & 1.609e-09 & 
3.99 \\
		\midrule
		\multicolumn{8}{c}{}\\[-3ex]
		\midrule
		{$\tau$} & {$h$} &
		{ $\| e^{u}  \|_{L^2(L^2)} $ } & EOC &
		{ $\| e^{v}  \|_{L^2(L^2)} $ } & EOC &
		{ $||| E\,  |||_{L^2} $ } & EOC \\
		\midrule
		$\tau_0/2^0$ & $h_0$ & 1.634e-04 & {--} & 1.232e-03 & {--} & 1.441e-03 & 
{--} \\
		$\tau_0/2^1$ & $h_0$ & 1.070e-05 & 3.93 & 7.864e-05 & 3.97 & 9.269e-05 & 
3.96 \\
		$\tau_0/2^2$ & $h_0$ & 6.765e-07 & 3.98 & 4.943e-06 & 3.99 & 5.836e-06 & 
3.99 \\
		$\tau_0/2^3$ & $h_0$ & 4.240e-08 & 4.00 & 3.094e-07 & 4.00 & 3.654e-07 & 
4.00 \\
		$\tau_0/2^4$ & $h_0$ & 2.652e-09 & 4.00 & 1.934e-08 & 4.00 & 2.285e-08 & 
4.00 \\
		$\tau_0/2^5$ & $h_0$ & 1.659e-10 & 4.00 & 1.212e-09 & 4.00 & 1.433e-09 & 
3.99 \\
		\bottomrule
	\end{tabular}
	\caption{Calculated errors for GCC$^1(3)$ with solution \eqref{eq:conv_test_1}.}
	\label{tab:conv_1}
\end{table}

In our second numerical experiment we study the space-time convergence behavior of a 
solution satisfying non-homogeneous Dirichlet boundary conditions,
\begin{equation}
\label{eq:u2}
u_2(\mat{x},t) = \sin(2 \cdot \pi \cdot t + x_1) \cdot \sin(2 \cdot \pi \cdot t \cdot 
x_2)
\end{equation} 
on $\Omega \times I = (0,1)^2 \times [0,1]$. We choose a bicubic discretization in 
$V_h^3$ (cf.\ \eqref{Eq:DefVh}) of the space variable. We refine the space-time mesh by 
halving both step sizes in each refinement step. Table~\ref{tab:conv_2} shows the computed 
errors and experimental orders of convergence for this example. In all measured norms, 
optimal rates in space and time (cf.\ Thm.~\ref{th:error}) are confirmed. This underlines 
the correct treatment of the prescribed non-homogeneous Dirichlet boundary conditions.
\begin{table}[!htb]
	\centering
	\small
	\begin{tabular}{c@{\,\,\,\,}c  c@{\,}c  c@{\,}c  c@{\,}c}
		\toprule
		{$\tau$} & {$h$} &
		{ $\| e^{u}  \|_{L^\infty(L^2)} $ } & EOC &
		{ $\| e^{v}  \|_{L^\infty(L^2)} $ } & EOC &
		{ $||| E\,  |||_{L^\infty} $ } & EOC \\
		\midrule
		$\tau_0/2^0$ & $h_0/2^0$ & 3.486e-03 & {--} & 3.602e-02 & {--} & 
5.013e-02 & {--} \\
		$\tau_0/2^1$ & $h_0/2^1$ & 2.329e-04 & 3.90 & 2.392e-03 & 3.90 & 
3.338e-03 & 3.92 \\
		$\tau_0/2^2$ & $h_0/2^2$ & 1.483e-05 & 3.97 & 1.527e-04 & 3.97 & 
2.128e-04 & 3.98 \\
		$\tau_0/2^3$ & $h_0/2^3$ & 9.320e-07 & 3.99 & 9.609e-06 & 3.99 & 
1.338e-05 & 3.99 \\
		$\tau_0/2^4$ & $h_0/2^4$ & 5.837e-08 & 4.00 & 6.022e-07 & 4.00 & 
8.383e-08 & 4.00 \\
		$\tau_0/2^5$ & $h_0/2^5$ & 3.649e-09 & 4.00 & 3.767e-08 & 4.00 & 
5.243e-08 & 4.00 \\
		\midrule
		\multicolumn{8}{c}{}\\[-3ex]
		\midrule
		{$\tau$} & {$h$} &
		{ $\| e^{u}  \|_{L^2(L^2)} $ } & EOC &
		{ $\| e^{v}  \|_{L^2(L^2)} $ } & EOC &
		{ $||| E\,  |||_{L^2} $ } & EOC \\
		\midrule
		$\tau_0/2^0$ & $h_0/2^0$ & 2.700e-03 & {--} & 2.568e-02 & {--} & 
3.458e-02 & {--} \\
		$\tau_0/2^1$ & $h_0/2^1$ & 1.771e-04 & 3.93 & 1.689e-03 & 3.93 & 
2.278e-03 & 3.92 \\
		$\tau_0/2^2$ & $h_0/2^2$ & 1.120e-05 & 3.98 & 1.070e-04 & 3.98 & 
1.444e-04 & 3.98 \\
		$\tau_0/2^3$ & $h_0/2^3$ & 7.020e-07 & 4.00 & 6.713e-06 & 3.99 & 
9.061e-06 & 3.99 \\
		$\tau_0/2^4$ & $h_0/2^4$ & 4.391e-08 & 4.00 & 4.199e-07 & 4.00 & 
5.669e-07 & 4.00 \\
		$\tau_0/2^5$ & $h_0/2^5$ & 2.744e-09 & 4.00 & 2.624e-08 & 4.00 & 
3.543e-08 & 4.00 \\
		\bottomrule
	\end{tabular}
	\caption{Calculated errors for GCC$^1(3)$ with solution \eqref{eq:u2}.}
	\label{tab:conv_2}
\end{table}

\subsection{Test case of structural health monitoring}
\label{sec:Challenging_Example}

Next, we consider a test problem that is based on \cite{Bangerth2010} and 
related to typical problems of structural health monitoring by ultrasonic waves (cf.\ \cref{fig:elastic_wave}). We aim to compare the GCC$^1(3)$ approach with a standard 
continuous in time Galerkin--Petrov approach cGP(2) of piecewise quadratic polynomials in 
time; cf.\ \cite{Bause2018,Hussain2011} for details. The cGP(2) scheme has superconvergence 
properties in the discrete time nodes $t_n$ for $n=1,\ldots, N$ as shown in 
\cite{Bause2018}. Thus,the errors  $\max_{n=1,\ldots, N} \|e^u(t_n)\|$ and 
$\max_{n=1,\ldots, N} \|e^v(t_n)\|$ for the GCC$^1(3)$ and the cGP(2) scheme admit the 
same fourth order rate of convergence in time and, thus, are comparable with respect to 
accuracy.

\begin{figure}[!htb]
\centering
\subcaptionbox{Test setting (according to \cite{Bangerth2010}). 	
\label{fig:Challenging_Plate}}
{\includegraphics[width=0.45\textwidth,keepaspectratio]
{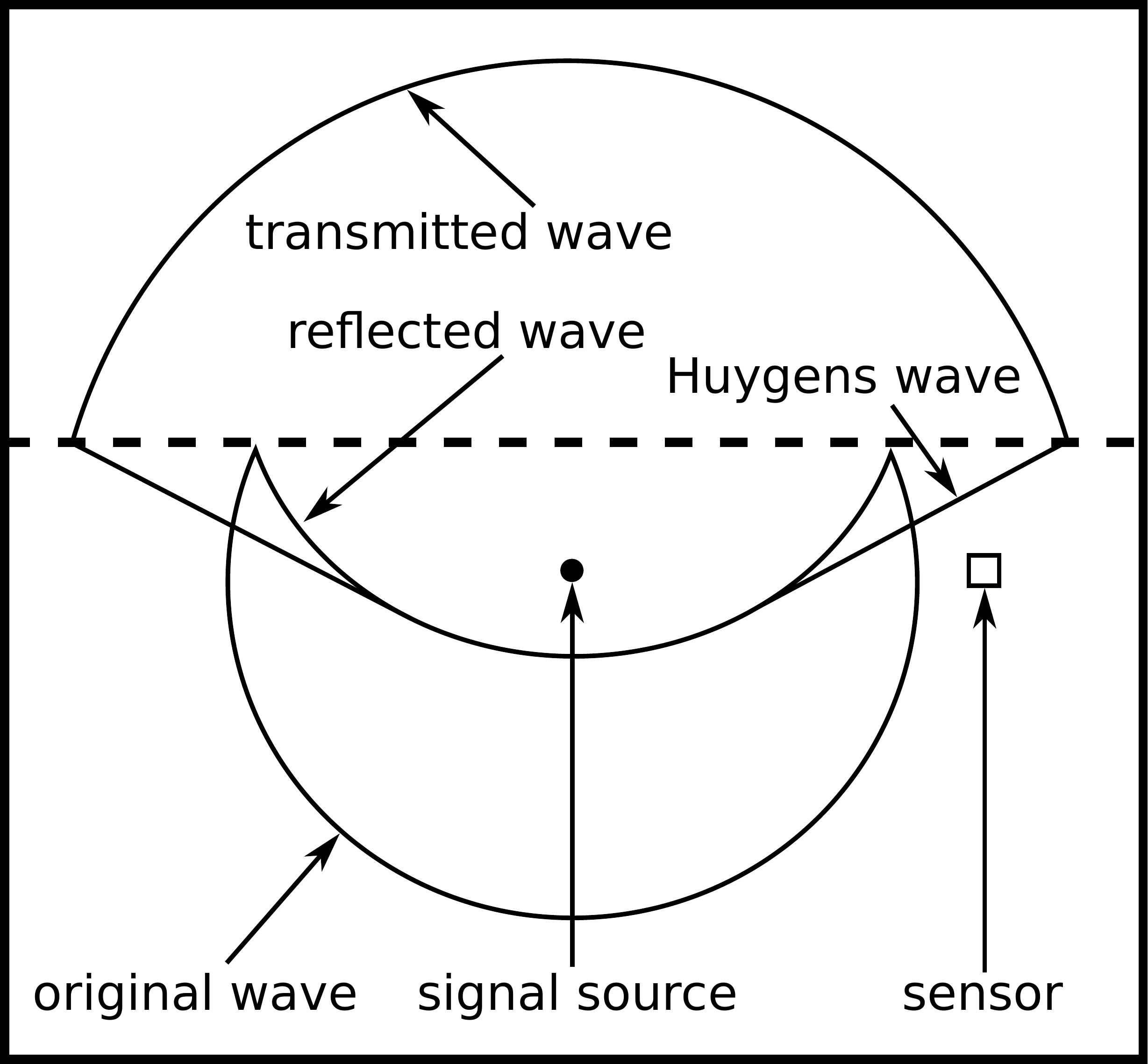}}
\subcaptionbox{Solution at $t=0.5$ with spatial mesh. \label{fig:Sol_t_0_5}}
[0.53\textwidth]
{\includegraphics[width=0.45\textwidth,keepaspectratio]
{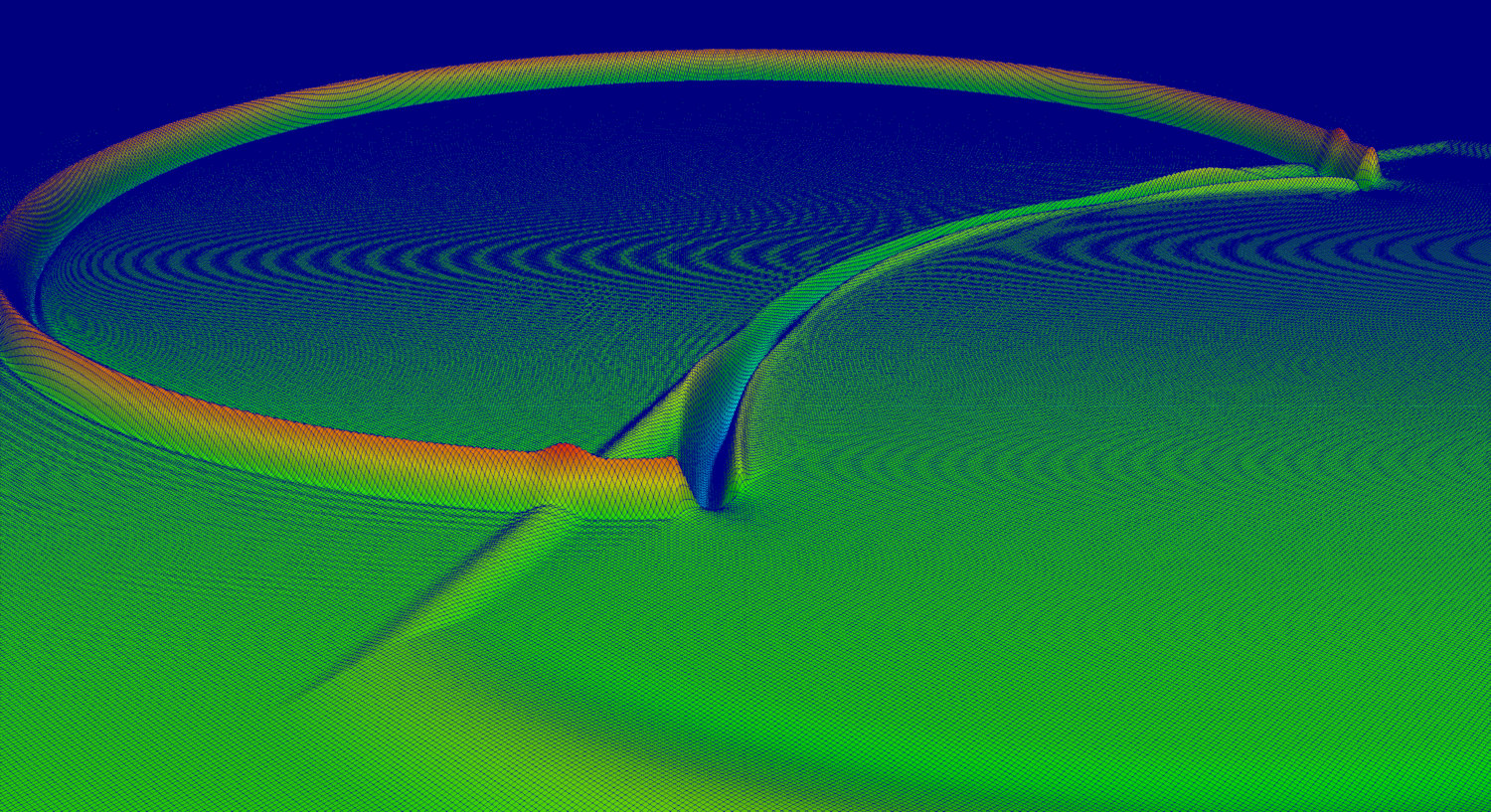}}
\caption{Test case of structural health monitoring}\label{fig:setup_preview}
\end{figure}

The test setting is sketched in \cref{fig:Challenging_Plate}. We consider $\Omega 
\times I = (-1,1)^2 \times (0,1)$, let $f = 0$ and, for simplicity, prescribe homogeneous 
Dirichlet boundary conditions such that $u^D = 0$. For the initial value we prescribe 
a regularized Dirac impulse by 
\begin{equation}
\begin{aligned}
u_0(\mat{x}) &= e^{-|\mat{x}_s|^2}
\big(1 - | \mat{x}_s |^2 \big) \Theta \big( 1 - | \mat{x}_s | \big),
&
\mat{x}_s &= 100 \mat{x}\,,
\end{aligned}
\end{equation}
where $\Theta$ is the Heaviside function. The coefficient function $c(\mat{x})$, mimicing 
a material parameter, has a jump discontinuity and is given by $c(\mat x) = 1$ for $x_2 
< 0.2$ and $c(\mat x) = 9$ for $x_2\geq 0.2$. Further we put $v_0=0$ for the second 
initial value. Finally, we define the control region $\Omega_c = (0.75 - h_c, 0.75 + h_c) 
\times (-h_c, h_c)$ where we calculate the signal arrival, at a sensor position for 
instance, in terms of 
\begin{equation}
\label{eq:control_quantity}
u_c(t) = \int_{\Omega_c} u_{\tau,h} (\mat{x}, t) \d \mat{x}.
\end{equation}

We choose a spatial mesh of \num{65536} cells and $\mathbb Q_7$ elements; cf.\ \cref{fig:Sol_t_0_5}. This leads to more than \num{3.2E6} degrees of freedom in space in 
each time step for each of the solution vectors. For each computation of the control 
quantity \eqref{eq:control_quantity}, with $t\in (0,1]$, we use a constant time step size 
$\tau_n$ for all time steps and compare the computation with the initially chosen 
reference time step size of $\tau_0 = \num{2e-5}$.

\begin{figure}[h!tbp]
\centering
\subcaptionbox{Control quantity \eqref{eq:control_quantity} for cGP(2) ({method $C^0$}) 
with different time step sizes and reference solution GCC$^1(3)$ ({method $C^1$}).}
{\fbox{\includegraphics[width=0.95\textwidth,height=1\textheight,keepaspectratio]
{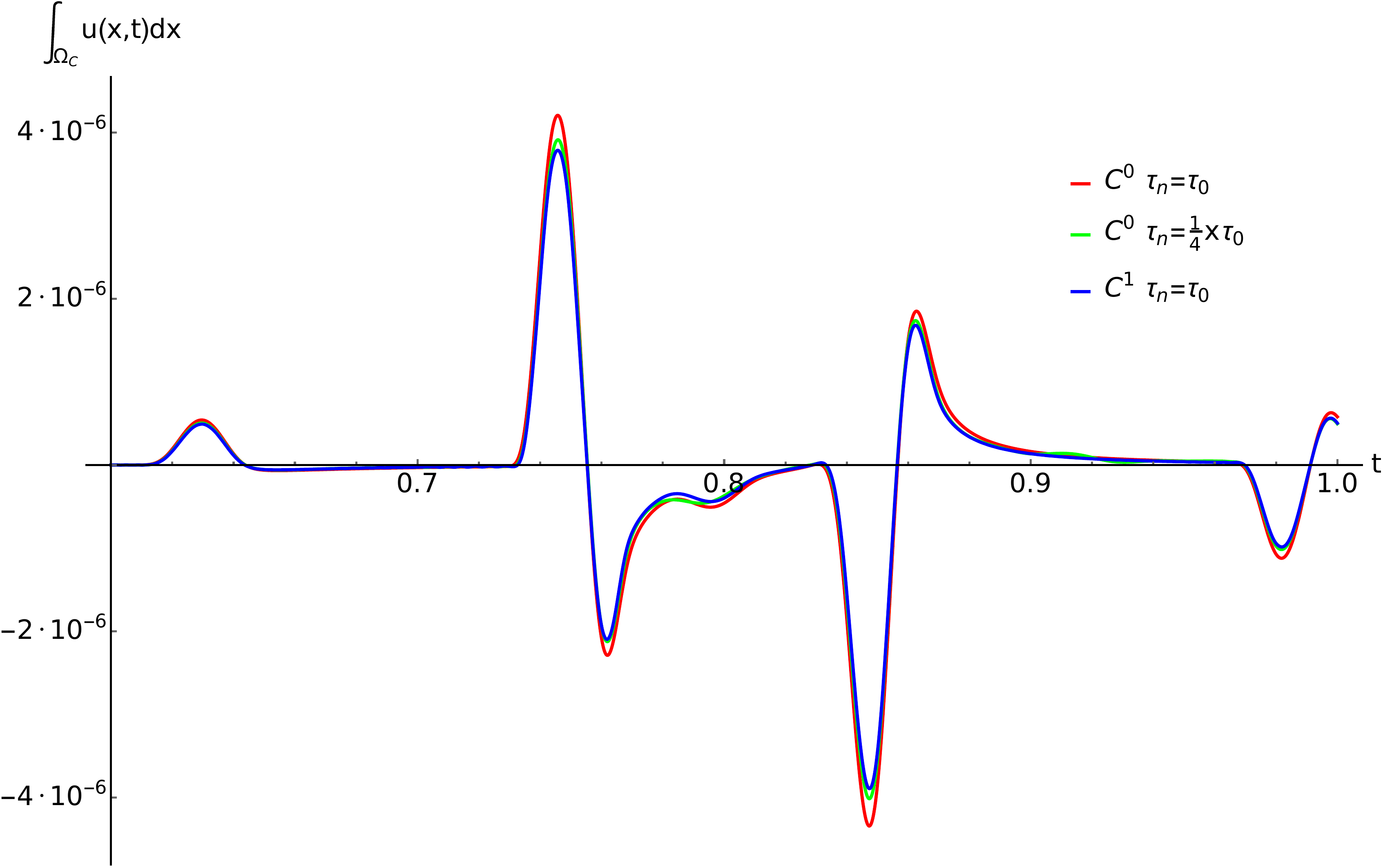}}}
\vspace{2ex}
\subcaptionbox{Control quantity \eqref{eq:control_quantity} for GCC$^1(3)$ ({method 
$C^1$}) with different time step sizes.}
{\fbox{\includegraphics[width=0.95\textwidth,height=1\textheight,keepaspectratio]
{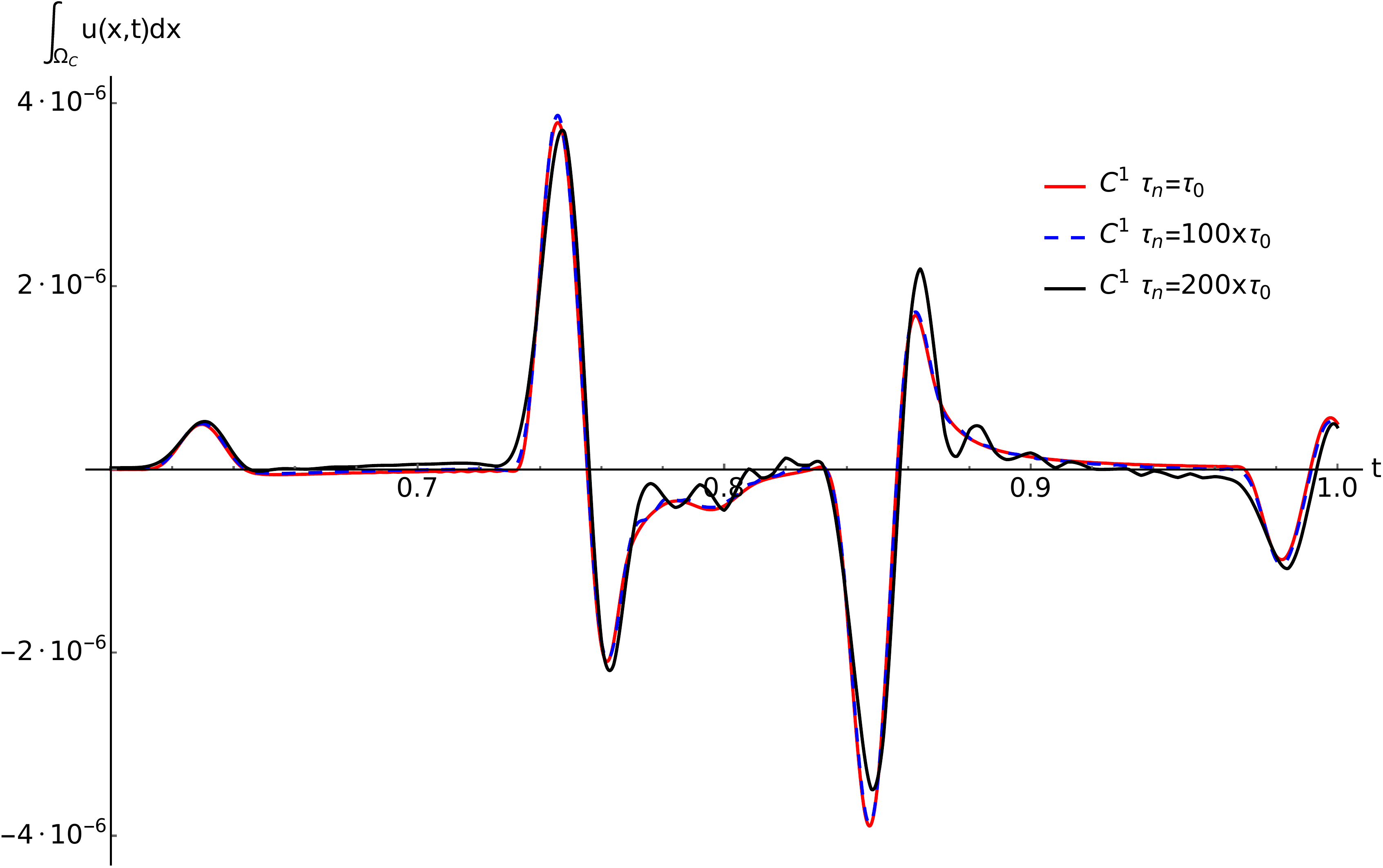}}}
\caption{Control quantity \eqref{eq:control_quantity} for GCC$^1(3)$ ({method 
$C^1$}) and cGP(2) ({method $C^0$}) for different time step sizes.}\label{fig:res_all}
\end{figure}

\Cref{fig:res_all} shows the signal arrival and control quantity 
\eqref{eq:control_quantity} over $t \in (0.6, 1)$ with different choices of the time 
step sizes for the Galerkin--collocation scheme GCC$^1(3)$ and the standard 
Galerkin--Petrov approach cGP(2) (cf.\ \cite{Bause2018,Hussain2011}) of a continuous in 
time approximation. For the cGP(2) approach, very small time step sizes are required to avoid over- and undershoots in the control quantity $u_c(t)$. For the GCC$^1(3)$ approach with $C^1$ regularity in time, much larger time steps, approximately 100 times larger, can be applied without loss of accuracy compared to the fully converged reference solution given by GCC$^1()3$ with step size $\tau_0$. This clearly shows the  superiority of the Galerkin--collocation scheme GCC$^1(3)$.

In Table~\ref{tab:runtime_comparison} the computational costs are summarized, where $r_1$ 
is the runtime for solving the condensed system by the approach of \cref{sec:method_1} and $r_2$ is the runtime for solving the block system by the 
approach of \cref{sec:method_2}. For the cGP(2) approach, only the first of 
the either iterative solver techniques was implemented. Recalling 
from \cref{fig:res_all} that GCC$^1(3)$ with $\tau_n = 100 \times \tau_0$ leads to 
the fully converged solution whereas cGP(2) with $\tau_n = \tau_0$ already shows over- 
and undershoots, a strong superiority of GCC$^1(3)$ over cGP(2) is observed in in 
Table~\ref{tab:runtime_comparison}. For both solver, a reduction in the wall clock 
time by a factor of about 25 is shown.  
\begin{table}[h!tb]
	\centering
	\small
	\begin{tabular}{cccrSS[table-number-alignment=center]}
		\toprule
		DoF (space) 	& cores	& method 	& \multicolumn{1}{c}{$\tau_n$}	& 
{$r_1$[h]} & {$r_2$[h]}\\
		\midrule
		\num{3.2E6}	& 224	& $C^0$	& $\num{0.25} \times \tau_0$ & 219.3 & 
{-} \\
		&		& $C^0$	& $\tau_0$ & 40.0 & {-} \\
		\addlinespace
		&		& $C^1$	& $\tau_0$ & 46.6 & 25.3 \\
		&		& $C^1$	& $2 \times \tau_0$ & 33.1 & 19.4 \\
		&		& $C^1$	& $25 \times \tau_0$ & 4.5 & 2.3 \\
		&		& $C^1$	& $35 \times \tau_0$ & 3.7 & 2.2 \\
		&		& $C^1$	& $50 \times \tau_0$ & 2.9 & 1.6\\
		&		& $C^1$	& $100 \times \tau_0$ & 1.7 & 0.9\\
		&		& $C^1$	& $200 \times \tau_0$ & 1.1 & 0.7 \\
		\addlinespace
		\num{4.2E6}	& 336	& $C^1$	& $50 \times \tau_0$ & 3.3 & 1.7\\
		\bottomrule
	\end{tabular}
	\caption{Runtime (wall clock time) for GCC$^1(3)$ ({method 
$C^1$}) and cGP(2) ({method $C^0$}) for different time step sizes and solvers 
of \cref{sec:method_1} ($r_1$) and \ref{sec:method_2} ($r_2$) .}
	\label{tab:runtime_comparison}
\end{table}

\section{\texorpdfstring{Galerkin--collocation GCC$\boldsymbol{{}^2(5)}$}
{Galerkin--collocation GCC2(5)}}
\label{sec:C2_solution}

Here, we briefly derive the algebraic form of the Galerkin--collocation scheme GCC$^2(k)$ of Def.~\ref{Def:GCC} with fully discrete solutions $(u_{\tau,h}{}_{|I_n},v_{\tau,h}{}_{|I_n})\in (X_\tau^k (V_h))^2$ such that $(u_{\tau,h},v_{\tau,h})\in (C^2(\overline I;V_h))^2$. For brevity, we restrict ourselves to the lowest polynomial degree in time $k=5$ that is possible to get $C^2$-regularity. The convergence properties are then demonstrated numerically. 

\subsection{Fully discrete system}

We follow the lines of \cref{sec:Deriving_C1_System} and use the notation introduced there. The six basis function of $\mathbb P_5 (\hat I;\R)$ on the reference interval $\hat I$ are defined by the conditions
\begin{align}
\hat{\xi_{i}}^{(l)}(j)
&=
\delta_{i-2*j-l, j}
&&
\forall \,
i \in \{0, \cdots, 5\}
\quad \land \quad
j \in \{0, 1\}\,,\quad
l \in \{0, 1, 2\}\,,
\end{align}
where $\delta_{i,j}$ denotes the usual Kronecker symbol. This gives us
\begin{equation}
\label{eq:basis_C2}
\mbox{}\hspace*{-3cm}
\begin{array}{r@{\;}c@{\;}lr@{\;}c@{\;}lr@{\;}c@{\;}l}
\hat{\xi_{0}} &= & -6 t^5+15 t^4-10 t^3+1\,,
& \hat{\xi_{1}} &= & -3 t^5+8 t^4-6 t^3+t\,,
& \hat{\xi_{2}} &= & -\frac{1}{2} t^5 +\frac{3}{2} t^4 - \frac{3}{2} t^3 + \frac{1}{2} t^2\,,\\[2ex]
\hat{\xi_{3}} &= & 6 t^5-15 t^4+10 t^3\,,
& \hat{\xi_{4}} &= & -3 t^5+7 t^4-4 t^3\,,
& \hat{\xi_{5}} &= & \frac{1}{2}t^5 - t^4 + \frac{1}{2} t^3\,.
\end{array}\hspace*{-3cm}\mbox{}
\end{equation}
For this basis of $\mathbb P_5(\hat I;\R)$, the discrete variational conditions \eqref{Eq:SemiDisLocalcGPC_9}, \eqref{Eq:SemiDisLocalcGPC_10} then read as 
\begin{align}
\label{eq:time_scheme_C2_a1}
& \mat{M} \biggl(
-\mat{u}_{n,0}^{0} + \mat{u}_{n,3}^{0}
\biggr)
-	\tau_n \mat{M} \biggl(
\frac{1}{2} \mat{v}_{n,0}^{0} + \frac{1}{10} \mat{v}_{n,1}^{0}
+ \frac{1}{120} \mat{v}_{n,2}^{0} + \frac{1}{2} \mat{v}_{n,3}^{0}
- \frac{1}{10} \mat{v}_{n,4}^{0} + \frac{1}{120} \mat{v}_{n,5}^{0}
\biggr)
= \mat 0\,,
\\[1ex]
&\begin{multlined}
\label{eq:time_scheme_C2_a2}
\mat{M} \biggl(
-\mat{v}_{n,0}^{0} + \mat{v}_{n,3}^{0}
\biggr)
+
\tau_n \mat{A}
\biggl(
\frac{1}{2} \mat{u}_{n,0}^{n,0} + \frac{1}{10} \mat{u}_{n,1}^{0}
+ \frac{1}{120} \mat{u}_{n,2}^{0} + \frac{1}{2} \mat{u}_{n,3}^{0}
- \frac{1}{10} \mat{u}_{n,4}^{0} + \frac{1}{120} \mat{u}_{n,5}^{0}
\biggr)
=
\\
\tau_n \mat{M}
\biggl( \frac{1}{2} \mat{f}_{n,0}^{} + \frac{1}{10} \mat{f}_{n,1}^{}
+ \frac{1}{120} \mat{f}_{n,2}^{} + \frac{1}{2} \mat{f}_{n,3}
- \frac{1}{10} \mat{f}_{n,4}^{} + \frac{1}{120} \mat{f}_{n,5}
\biggr)
-
\mat{M} \biggl(
-\mat{v}_{n,0}^{D} + \mat{v}_{n,3}^{D}
\biggr)
\\
-
\tau_n \mat{A}
\biggl(
\frac{1}{2} \mat{u}_{n,0}^{D} + \frac{1}{10} \mat{u}_{n,1}^{D}
+ \frac{1}{120} \mat{u}_{n,2}^{D} + \frac{1}{2} \mat{u}_{n,3}^{D}
- \frac{1}{10} \mat{u}_{n,4}^{D} + \frac{1}{120} \mat{u}_{n,5}^{D}
\biggr)\,.
\end{multlined}
\end{align}
In the basis, the first collocation conditions \eqref{Eq:SemiDisLocalcGPC_6} yield for $\mat{w}_{n,i}^0 \in \{\mat{u}_{n,i}^0,  \mat{v}_{n,i}^0\}$ that 
\begin{align}
\label{eq:constraint_left_C2_a}
\mat{w}_{n,0}^{0} &= \mat{w}_{n-1,3}^{0}\,,
&
\mat{w}_{n,1}^{0} &= \mat{w}_{n-1,4}^{0}\,,
&
\mat{w}_{n,2}^{0} &= \mat{w}_{n-1,5}^{0}\,,
\end{align}
which reduces the number of unknown solution vectors by $6$ on each subinterval $I_n$. For the collocation conditions \eqref{Eq:SemiDisLocalcGPC_7}, \eqref{Eq:SemiDisLocalcGPC_8} at $t_n$ and $s=1$ we deduce that 
\begin{align}
\label{eq:C2_Col2a1}
&\mat{M}
\frac{1}{\tau_n} \mat{u}_{n,4}^{0}
-
\mat{M}
\mat{v}_{n,3}^{0}
= \mat 0\,,\quad 
\begin{split}
\label{eq:C2_Col2a2}
& \mat{M} \frac{1}{\tau_n} \mat{v}_{n,4}^{0}
+
\mat{A} \mat{u}_{n,3}^{0}
=
\mat{M}
\mat{f}_{n,3}^{}
-
\mat{M} \frac{1}{\tau_n} \mat{v}_{n,4}^{D}
-
\mat{A} \mat{u}_{n,3}^{D}\,.
\end{split}
\end{align}
Similarly, for $s=2$ the collocation conditions \eqref{Eq:SemiDisLocalcGPC_7}, \eqref{Eq:SemiDisLocalcGPC_8} at $t_n$ read as 
\begin{align}
\label{eq:C2_Col2b1}
& \mat{M}
\frac{1}{\tau_n} \mat{u}_{n,5}^{0}
-
\mat{M}
\mat{v}_{n,4}^{0}
=\mat 0\,, \quad 
\begin{split}
\label{eq:C2_Col2b2}
&\mat{M} \frac{1}{\tau_n} \mat{v}_{n,5}^{0}
+
\mat{A} \mat{u}_{n,4}^{0}
=
\mat{M}
\mat{f}_{n,4}^{}
-
\mat{M} \frac{1}{\tau_n} \mat{v}_{n,5}^{D}
-
\mat{A} \mat{u}_{n,4}^{D}.
\end{split}
\end{align}
Finally, we recover the previous conditions as the linear system $\mat S \mat x = \mat 
b$ for the vector of unknowns $\mat x =\left((\mat{u}_{n,3}^{0})^\top, 
(\mat{u}_{n,4}^{0})^\top, (\mat{v}_{n,5}^{0})^\top, 
\mat{u}_{n,5}^{0})^\top\right)^\top$ and with the system matrix $\mat S$ and right-hand side 
vector $\mat b$ given by
\begingroup
\renewcommand*{\arraystretch}{1.5}
\begin{align}
\label{eq:system_matrix_c2}
\mat{S} &=
\begin{pmatrix}
\mat{A} & \mat{0} & \mat{0} & \frac{1}{\tau_n^2}\mat{M} \\
\mat{0} & \mat{A} & \frac{1}{\tau_n}\mat{M} & \mat{0} \\
\mat{M} & -\frac{1}{2}\mat{M} & -\frac{\tau_n}{120} \mat{M} & \frac{1}{10}\mat{M} \\
\frac{\tau_n}{2}\mat{A} & \frac{1}{\tau_n}\mat{M} -\frac{\tau_n}{10} \mat{A} & \mat{0} & 
\frac{\tau_n}{120}\mat{A}
\end{pmatrix}.
&
\mat{b} &= 
\begin{pmatrix}
\mat{f}_{n,3} - \mat{A} \mat{u}_{n,3}^D - \frac{1}{\tau_n} \mat{M} \mat{v}_{n,4}^D
\\
\mat{f}_{n,4} - \mat{A} \mat{u}_{n,4}^D - \frac{1}{\tau_n} \mat{M} \mat{v}_{n,5}^D
\\
\mat{b}_{n,3}
\\
\mat{b}_{n,4}
\end{pmatrix},
\end{align}
\endgroup
with $\mat{b}_{n,3} =
\mat{M} \bigl(
\mat{u}_{n,0}^0 + \mat{u}_{n,0}^D - \mat{u}_{n,3}^D
+ \frac{\tau_n}{2} (\mat{v}_{n,0}^0 + \mat{v}_{n,0}^D)
+ \frac{\tau_n}{10} (\mat{v}_{n,1}^0 + \mat{v}_{n,1}^D)
+ \frac{\tau_n}{120} (\mat{v}_{n,2}^0 + \mat{v}_{n,2}^D)
+ \tau_n (\frac{1}{2} \mat{v}_{n,3}^D - \frac{1}{10} \mat{v}_{n,4}^D + \frac{1}{120} \mat{v}_{n,5}^D)
$
and
$
\mat{b}_{n,4} =
\mat{M} \bigl( \mat{v}_{n,0}^0 + \mat{v}_{n,0}^D - \mat{v}_{n,3}^D \bigr)
+ \tau_n \bigl( \frac{1}{2} \mat{f}_{n,3}  + \frac{1}{10} \mat{f}_{n,1} + \frac{1}{120} \mat{f}_{n,2} + \frac{1}{2} \mat{f}_{n,3} - \frac{1}{10} \mat{f}_{n,4} + \frac{1}{120} \mat{f}_{n,5}\bigr)
- 
\tau_n \mat{A} \bigl(
\frac{1}{2} (\mat{u}_{n,0}^0 + \mat{u}_{n,0}^D)2
+ \frac{1}{10} (\mat{u}_{n,1}^0 + \mat{u}_{n,1}^D)
+ \frac{1}{120} (\mat{u}_{n,2}^0 + \mat{u}_{n,2}^D)
+ \frac{1}{2} \mat{u}_{n,3}^D
- \frac{1}{10} \mat{u}_{n,4}^D
+ \frac{1}{120} \mat{u}_{n,5}^D
\bigr)
$.

\subsection{Iterative solver and convergence study}
\label{sec:NumTestC2}

To solve the linear system $\mat S \mat x = \mat b$ with $\mat{S}$ from \eqref{eq:system_matrix_c2}, we use block Gaussian elimination, as sketched in \cref{sec:method_1}, to find a reduced system $\mat S_r \mat{u}_{n,4}^{0} = \mat b_r$ for the essential unknown $\mat{u}_{n,4}^{0}$. All remaining unknown subvectors of $\mat x$ can be computed in post-processing steps. In explicit form, the condensed system reads as
\begin{equation}
\label{eq:recuced_c2_system}
\left(
14400 \mat{M} + 720 \tau_n^2 \mat{A} + 24 \tau_n^4 \mat{A}\mat{M}^{-^1}\mat{A} + \tau_n^6 
\mat{A}\mat{M}^{-^1}\mat{A}\mat{M}^{-^1}\mat{A}
\right)
\mat{u}_{n,4}^{0}
=
\mat{b}_r\,.
\end{equation}
For brevity, we omit the exact definition of $\mat b_r$ that can  be deduced easily  from 
\eqref{eq:system_matrix_c2}.

The matrix $\mat S_r$ is symmetric such that preconditioned conjugate gradient iterations are used for its solution. The preconditioner is constructed along the lines of \cref{sec:method_1}. The remainder part $\tau_n^6 \mat{A}\mat{M}^{-^1}\mat{A}\mat{M}^{-^1}\mat{A}$ is still ignored in the construction of the preconditioner. Even though the remainder is weighted by the small factor $\tau_n^6$, numerical experiments indicate that this scaling is not sufficient to balance its impact on the interation process. For the construction of an efficient preconditioning technique for $\mat S_r$ of GCC$^2(5)$ further improvements are still necessary.   

To illustrate the convergence behavior and performance of the GCC$^2(5)$ Galerkin--collocation approach, we present in Table~\ref{tab:conv_C2_1} our numerical results for the test problem \eqref{eq:conv_test_1}. The expected convergence of sixth order in time is nicely observed in all norms.
\begin{table}[!htb]
	\centering
	\small
	\begin{tabular}{c@{\,\,\,\,}c  c@{\,}c  c@{\,}c  c@{\,}c}
		\toprule
		{$\tau$} & {$h$} &
		{ $\| e^{u}  \|_{L^\infty(L^2)} $ } & EOC &
		{ $\| e^{v}  \|_{L^\infty(L^2)} $ } & EOC &
		{ $||| E\,  |||_{L^\infty} $ } & EOC \\
		\midrule
		$\tau_0/2^0$ & $h_0$ & 8.748e-06 & {--} & 4.355e-05 & {--} & 4.985e-05 & 
{--} \\
		$\tau_0/2^1$ & $h_0$ & 1.370e-07 & 6.00 & 7.404e-07 & 5.88 & 8.043e-07 & 
5.95 \\
		$\tau_0/2^2$ & $h_0$ & 2.165e-09 & 5.98 & 1.202e-08 & 5.95 & 1.266e-08 & 
5.99 \\
		$\tau_0/2^3$ & $h_0$ & 3.388e-11 & 6.00 & 1.883e-10 & 6.00 & 1.980e-10 & 
6.00 \\
		$\tau_0/2^4$ & $h_0$ & 5.301e-13 & 6.00 & 2.940e-12 & 6.00 & 3.093e-12 & 
6.00 \\
		%
		\midrule
		\multicolumn{8}{c}{}\\[-3ex]
		\midrule
		{$\tau$} & {$h$} &
		{ $\| e^{u}  \|_{L^2(L^2)} $ } & EOC &
		{ $\| e^{v}  \|_{L^2(L^2)} $ } & EOC &
		{ $||| E\,  |||_{L^2} $ } & EOC \\
		\midrule
		$\tau_0/2^0$ & $h_0$ & 4.022e-06 & {--} & 2.996e-05 & {--} & 3.502e-05 & 
{--} \\
		$\tau_0/2^1$ & $h_0$ & 6.353e-08 & 5.98 & 4.808e-07 & 5.96 & 5.599e-07 & 
5.97 \\
		$\tau_0/2^2$ & $h_0$ & 9.957e-10 & 6.00 & 7.565e-09 & 5.99 & 8.800e-09 & 
5.99 \\
		$\tau_0/2^3$ & $h_0$ & 1.557e-11 & 6.00 & 1.184e-10 & 6.00 & 1.377e-10 & 
6.00 \\
		$\tau_0/2^4$ & $h_0$ & 2.431e-13 & 6.00 & 1.849e-12 & 6.00 & 2.151e-12 & 
6.00 \\
		%
		\bottomrule
	\end{tabular}
	\caption{Calculated errors for GCC$^2(5)$ with solution \eqref{eq:conv_test_1}.}
	\label{tab:conv_C2_1}
\end{table}

\vspace*{-6ex}


\end{document}